\documentclass[10pt]{amsart}
\usepackage{amssymb}
\usepackage{amsmath}
\usepackage{amsthm}
\usepackage{graphicx}
\usepackage{subfigure}
\usepackage{url}
\usepackage{hyperref}
\usepackage{color}
\usepackage{fullpage}
\usepackage{cite}

\numberwithin{equation}{section}

\newcommand{\abs}[1]{\left| #1\right|}
\newcommand{\norm}[1]{\left\|#1\right\|}
\newcommand{\paren}[1]{\left(#1\right)}
\newcommand{\bracket}[1]{\left[#1\right]}
\newcommand{\set}[1]{\left\{#1\right\}}

\newcommand{\TV}{\mathrm{TV}}

\newcommand{\bigo}{\mathrm{O}}

\newcommand{\mean}[1]{\left\langle#1\right\rangle}

\newcommand{\refe}{\mathrm{ref}}
\newcommand{\simu}{\mathrm{sim}}
\newcommand{\phase}{\mathrm{phase}}

\newcommand{\tcorr}{t_{\mathrm{corr}}}
\newcommand{\tphase}{t_{\mathrm{phase}}}

\newcommand{\tpar}{t_{\mathrm{par}}}

\newcommand{\obs}{\mathcal{O}}

\newcommand{\tol}{{\rm TOL}}

\newcommand{\E}{\mathbb{E}}
\newcommand{\prob}{\mathbb{P}}

\newcommand{\R}{\mathbb{R}}

\newcommand{\LJ}{{\rm LJ}_7^{{\rm 2D}}}
\newcommand{\Var}{{\rm Var}}

\begin{document}
\title[Generalized ParRep]{A Generalized Parallel Replica Dynamics}

\author[Binder]{Andrew Binder} \author[Leli\`evre]{Tony Leli\`evre}
\author[Simpson]{Gideon Simpson}

\date{\today}

\maketitle

\begin{abstract}
  Metastability is a common obstacle to performing long molecular
  dynamics simulations. Many numerical methods have been proposed to
  overcome it.  One method is parallel replica dynamics, which relies
  on the rapid convergence of the underlying stochastic process to a
  quasi-stationary distribution.  Two requirements for applying
  parallel replica dynamics are knowledge of the time scale on which
  the process converges to the quasi-stationary distribution and a
  mechanism for generating samples from this distribution.  By
  combining a Fleming-Viot particle system with convergence
  diagnostics to simultaneously identify when the process converges
  while also generating samples, we can address both points.  This
  variation on the algorithm is illustrated with various numerical
  examples, including those with entropic barriers and the 2D
  Lennard-Jones cluster of seven atoms.
\end{abstract}


\section{Introduction}
\label{s:intro}

An outstanding obstacle for many problems modeled by {\it in situ}
molecular dynamics (MD) is the vast separation between the
characteristic time for atomic vibrations ($10^{-15}$ s), and the
characteristic time for macroscopic phenomena ($10^{-9}$ -- $10^{-3}$
s).  At the heart of this scale separation is the presence of {\it
  metastable regions} in the configuration space of the problem.
Examples of metastable configurations include the defect arrangement
in a crystal or the conformation of a protein.  Such metastability may
be due to either the energetic barriers of a potential energy driving
the problem or to the entropic barriers arising from steric
constraints. In the first case (energetic barriers), metastability is
due to the system needing to pass through a saddle point which is
higher in energy than the local minima to get from one metastable
region to another. In the second case (entropic barriers),
metastability is due to the system having to find a way through a
narrow (but not necessarily high energy) corridor to go from one large
region to another (see Section \ref{s:ent2d} below for an example with
entropic barriers)

Motivated by the challenge of this time-scale separation, A.F. Voter
proposed several methods to conquer metastability in the 1990s:
Parallel Replica Dynamics (ParRep), Temperature Accelerated Dynamics
(TAD) and Hyperdynamics (Hyper),
\cite{perez2009accelerated,Sorensen:2000p13728,Voter:1997p12684,Voter:1997p12731,Voter:1998p13729,Voter:2002p12678}.
These methods were derived using Transition State Theory and intuition
developed from kinetic Monte Carlo models, as the latter describes the
hopping dynamics between metastable regions. Indeed, the aim of all
these algorithms is to efficiently generate a realization of the
discrete-valued jump process amongst metastable regions. The main idea
is that the details of the dynamics within each metastable region are
not essential to our physical understanding. Rather, the goal should
be to get the correct statistics of the so-called {\em state-to-state
  dynamics}, corresponding to jumps amongst the metastable
regions. This is nontrivial in general for two reasons: (i)~the
original dynamics projected onto the state-to-state dynamics are not
Markovian; (ii) the parameters (transition rates) of the underlying
state-to-state dynamics are unknown.

In recent mathematical studies of these approaches, it has been shown
that these three algorithms take advantage of {\it quasi-stationary
  distributions} (QSDs) associated with the metastable states,
see~\cite{Aristoff:2014ch,LeBris:2012et,Lelievre:fk,Simpson:2013cs}.
Crudely, the QSD corresponds to the distribution of the end points of
trajectories conditioned on persisting in the region of interest for a
very long time. {This mathematical formalization clarifies
  the fundamental assumptions under which the algorithms will be
  accurate and also broadens their applicability. Indeed, one aim of
  this paper is to propose a modification of the original ParRep
  algorithm that will allow for states to be defined by generic
  partitions of configuration space.}  This bears some resemblance to
{\it milestoning}, which also allows more general partitions of
configuration space, \cite{Kirmizialtin:2011ba,VandenEijnden:2008bn}.
Therefore, we will not refer to ``basins of attractions'' or
``metastable regions'', but rather simply to ``states''. The only
requirement is that these states define a partition of the
configuration space. The boundary at the interface between two states
is called the dividing surface.

Briefly (this is detailed in Section~\ref{sec:orig_parrep} below),
ParRep works by first allowing a single reference trajectory to
explore a state.  If the trajectory survives for sufficiently long,
its end point will agree, in law, with the aforementioned QSD.  One
thus introduces a {\it decorrelation time}, denoted $\tcorr$, as the
time at which the law of the reference process will have converged to
the QSD.  Provided the reference process survives in the state up till
$\tcorr$, it is replaced by an ensemble of $N$ independent and
identically distributed replicas, each with an initial condition drawn
from the QSD.  The first replica to escape is then followed into the
next state.  As the replicas evolve independently and only a first
escape is desired, they are readily simulated in parallel, providing
as much as a factor of $N$ speedup of the exit event. Thus, there are
two practical challenges to implementing ParRep:
\begin{itemize}
\item Identifying a $\tcorr$ at which the law of the reference process
  is close to the QSD.
\item Generating samples from the QSD from which to start the
  replicas.
\end{itemize}
In the original algorithm, $\tcorr$ is {\em a priori} chosen by the
user, as the states are defined so that an approximation of the time
required to get ``local equilibration within the state'' is available.
{Such a value can be estimated in the case of energetic
  barriers at sufficiently low temperature using harmonic transition
  state theory.  But in the case of entropic barriers, there is, in
  general, no simple estimate; see, however, \cite{Kum:2004iu} for an
  example where ParRep was applied to an entropic barrier as an
  estimate of $\tcorr$ was available.  In the original algorithm,} the
sampling of the QSD is done using a rejection algorithm, which will be
inefficient if the state does not correspond to a metastable region
for the original dynamics. Indeed, this will degenerate in the long
time limit, as the trajectories will always exit.  In this work, we
propose an algorithm addressing both points, based on two ingredients:
\begin{itemize}
\item The use of a branching and interacting particle system called
  the Fleming-Viot particle process to simulate the law of the process
  conditioned on persisting in a state, and to sample the QSD in the
  longtime limit.
\item The use of Gelman-Rubin statistics in order to identify the
  correlation time, namely the convergence time to a stationary state
  for the Fleming-Viot particle process.
\end{itemize}
As we state below, this modified version of ParRep, presented below in
Section \ref{sec:modParRep}, relies on assumptions (see (A1) and (A2)
below) which would require more involved analysis to fully justify.
We demonstrate below in a collection of numerical experiments that
this modified algorithm gives results consistent with direct
simulations. We observe speedup factors up to ten in our test problems
where $N=100$ replicas were used. The aim of this paper is to present
new algorithmic developments, and not to explore the mathematical
foundations underpinning these ideas.

Though we focus on the ParRep algorithm, since it is the most natural
setting for introducing the Fleming-Viot particle process, identifying
the convergence to the QSD is also relevant to other problems,
including the two other accelerated dynamics algorithms: Hyper and
TAD. See~\cite{Aristoff:2014ch,Lelievre:fk} for the relevant
discussions.

Our paper is organized as follows.  In Section \ref{sec:orParRep}, we
introduce the dynamics of interest, review some properties of the QSD
and recall the original ParRep algorithm.  In
Section~\ref{sec:modParRep}, we then present the Fleming-Viot particle
process and the Gelman-Rubin statistics, which are needed to build the
modified ParRep algorithm we propose. Finally, in
Section~\ref{s:offline_examples} we show the effectiveness and caveats
of convergence diagnostics, before exploring the efficiency and
accuracy of the modified ParRep algorithm on various test cases in
Section~\ref{s:parrep_examples}.

\subsection{Acknowledgments}

A.B. was supported by a US Department of Defense NDSEG fellowship.
G.S. was supported in part by the US Department of Energy Award
DE-SC0002085 and the US National Science Foundation PIRE Grant
OISE-0967140.  T.L. acknowledges funding from the European Research
Council under the European Union's Seventh Framework Programme
(FP7/2007-2013) grant agreement no. 614492.

\noindent The authors would also like to thank C. Le Bris, M. Luskin,
D. Perez, and A.F. Voter for comments and suggestions throughout the
development of this work.  The authors would also like to thank the
referees for their helpful remarks.

\section{The Original ParRep and Quasi-Stationary Distributions}
\label{sec:orParRep}

\subsection{Overdamped Langevin Dynamics}
\label{s:odlang_dynamics}

We consider the case of the overdamped Langevin equation
\begin{equation}
  \label{e:odlang}
  dX_t = - \nabla V(X_t) dt + \sqrt{2\beta^{-1}} dW_t.
\end{equation}
Here, the stochastic process $(X_t)_{t \ge 0}$ takes values in $\R^d$,
$\beta$ is the inverse temperature, $V(x)$ is the driving potential
and $W_t$ a standard $d$-dimensional Brownian motion.  In all that
follows, we focus, for simplicity, on \eqref{e:odlang}.  However, the
algorithm we propose equally applies to the phase-space Langevin
dynamics which are also of interest. As mentioned in the introduction,
for typical potentials, the stochastic process $(X_t)_{t \ge 0}$
satisfying~\eqref{e:odlang} is metastable.  Much of its trajectory is
confined to particular regions of $\R^d$, occasionally hopping amongst
them.


Let $\Omega\subset \R^d$ denote the region of interest (namely the
state), and define
\begin{equation}
  \label{e:firstexit}
  T = \inf\set{t\geq 0\mid X_t \notin \Omega}
\end{equation}
to be the {\it first exit time} from $\Omega$, where $X_0 = x\in
\Omega$.  The point on the boundary, $X_T \in \partial \Omega$, is the
{\it first hitting point}. The aim of accelerated dynamics algorithms
(and ParRep in particular) is to efficiently sample $(T,X_T)$ from the
exit distribution.

\subsection{Quasi-stationary Distributions}
\label{s:qsd}

In order to present the original ParRep algorithm, it is helpful to be
familiar with quasi-stationary distributions (QSD). For more details
about quasi-stationary distributions, we refer the reader to, for
example,
~\cite{Cattiaux:2009um,Cattiaux:2010fx,Collet:1995ts,Martinez:1994vn,Martinez:2004we,Steinsaltz:2007cy,collet13:_quasi_station_distr}. Reference
\cite{LeBris:2012et} gives self-contained proofs of the results below.

Consider a smooth bounded open set $\Omega \subset \R^d$, that
corresponds to {a state}. By definition, the quasi-stationary
distribution $\nu$, associated with the dynamics~\eqref{e:odlang} and
the state $\Omega$, is the probability distribution, with support on
$\Omega$, satisfying, for all (measurable) $A\subset \Omega$ and
$t\geq 0$,
\begin{equation}
  \label{e:qsd}
  \begin{split}
    \nu(A) = \frac{\int \prob^{x} \bracket{X_t\in A,
        T>t}\nu(dx)}{\int\prob^{x} \bracket{T>t}\nu(dx)} &=
    \frac{\prob^{\nu} \bracket{X_t\in A, T>t}}{\prob^{\nu}
      \bracket{T>t}}= \prob^{\nu}\bracket{X_t \in A\mid T>t}.
  \end{split}
\end{equation}
Here and in the following, we indicate by a superscript the initial
condition for the stochastic process: $\prob^{x}$ indicates that
$X_0=x$ and $\prob^{\nu}$ indicates that $X_0$ is distributed
according to $\nu$.  In our setting, it can be shown that $\nu$ exists
and is unique.

The QSD enjoys three properties. First, it is related to an elliptic
eigenvalue problem. Let $L$ be the infinitesimal generator
of~\eqref{e:odlang}, defined by, for any smooth function $v:\R^d \to
\R$,
\begin{equation}
  \label{e:inf_gen}
  L v  =- \nabla V\cdot \nabla v + \beta^{-1} \Delta v.
\end{equation}
The operator $L$ is related to the stochastic process through the
following well-known result: if the function $u:\R^+ \times \Omega \to
\R$ satisfies the Kolmogorov equation:
\begin{equation}
  \label{e:kolmogorov}
  \left\{
    \begin{aligned}
      \partial_t u &= L u = - \nabla V\cdot \nabla u + \beta^{-1}
      \Delta
      u, \quad \text{for $t > 0$, $x \in  \Omega$}, \\
      u(t,x)&= f(x) \quad \text{for $t > 0$, $x \in \partial \Omega$},\\
      u(0,x)&=u_0(x) \quad \text{for $x \in\Omega$},
    \end{aligned}
  \right.
\end{equation}
then $u$ admits the probabilistic representation formula (Feynman-Kac
relation):
\begin{equation}
  \label{e:fk}
  u(t,x) = \E^x \bracket{u_0(X_t) 1_{T>t}} + \E^x\bracket{f(X_T)
    1_{T\leq t}}.
\end{equation}
Recall that $T$ defined by \eqref{e:firstexit}, is the first exit time
of $X_t$ from $\Omega$.  Provided $\Omega$ is bounded with
sufficiently smooth boundary and $V$ is smooth, $L$ has an infinite
set of Dirichlet eigenvalues and orthonormal eigenfunctions
\begin{equation}
  \label{e:spectrum}
  L \varphi_j = -\lambda_j \varphi_j, \quad \varphi_j|_{\partial \Omega}
  = 0, \quad j=1,2,\ldots
\end{equation}
{Here, the eigenfunctions are orthonormal with respect to the invariant
measure restricted to $\Omega$:
\begin{equation*}
\frac{\int_\Omega \varphi_j(x) \varphi_k(x) \exp(-\beta V(x))
  dx}
{\int_\Omega \exp(-\beta V(x))
  dx} = \delta_{jk}, \quad \forall j,k=1,2,\ldots
\end{equation*}}
The eigenfunction associated with the lowest eigenvalue is signed,
and, taking it to be positive, the QSD is
\begin{equation}
  \label{e:qsd_eig}
  \nu(dx) = \frac{\displaystyle \varphi_1(x) e^{-\beta V(x)}dx}{\displaystyle \int_\Omega \varphi_1(x) e^{-\beta V(x)}dx}.
\end{equation}
While this expression is explicit in terms of $\varphi_1$ and $V$,
since the problem is posed in $\R^d$ with $d$ large, it is not
practical to sample the QSD directly by first computing $\varphi_1$.

The second property associated to the QSD is that for all $t\geq 0$
and $A\subset \partial \Omega$:
\begin{equation}
  \label{e:qsd_indep}
  \begin{split}
    \prob^{\nu}\bracket{X_T\in A, T>t} &= \prob^{\nu}\bracket{
      T>t}\, \prob^{\nu}\bracket{X_T\in A} \\
    &=\paren{e^{-\lambda_1 t}}\, \paren{\int_A -\frac{1}{\beta
        \lambda_1} \frac{\nabla \varphi_1 e^{-\beta V} \cdot{\bf
          n}}{{\int_\Omega \varphi_1(x) e^{-\beta V(x)}dx}} d{S} },
  \end{split}
\end{equation}
where $dS$ is the surface Lebesgue measure on $\partial \Omega$ and
${\bf n}$ the unit outward normal vector to $\Omega$. Thus, the first
hitting point and first exit time are independent, and exit times are
exponentially distributed.  These two properties will be one of the
main arguments justifying ParRep. As explained
in~\cite{LeBris:2012et}, they are consequences of~\eqref{e:fk}
and~\eqref{e:qsd_eig}.

The third property of the QSD also plays an important role in
ParRep. Let us again consider $X_t$ satisfying~\eqref{e:odlang} with
$X_0 \sim \mu_0$ ($\mu_0$ with support in $\Omega$). Define the law of
$X_t$, conditioned on non-extinction as:
\begin{equation}
  \label{e:conditional_law}
  \mu_t(\bullet) =   \frac{\prob^{\mu_0}\bracket{X_t \in \bullet,
      T>t} }{\prob^{\mu_0}\bracket{ T>t} } =  \prob^{\mu_0}\bracket{X_t \in \bullet\mid T>t}.
\end{equation}
One can check that
\begin{equation}
  \mu_t(v_0) = \E^{\mu_0}\bracket{v_0(X_t)\mid T>t} = \frac{\int_\Omega v(x,t) \mu_0(dx)}{\int_\Omega \bar
    v(x,t) \mu_0(dx)},
\end{equation}
where $v$ solves \eqref{e:kolmogorov} with initial condition $v(0,x) =
v_0(x)$ and boundary conditions $v|_{\partial\Omega} = 0$ while $\bar
v$ solves \eqref{e:kolmogorov} with initial condition $\bar v(0,x) =
1$ and boundary conditions $\bar v|_{\partial\Omega} = 0$. Through
eigenfunction expansions of the form
\begin{equation}
  v(t,x) = \sum_{k=1}^\infty e^{-\lambda_k t} \varphi_k(x) \int v_0(y) \varphi_k(y)
  \frac{e^{-\beta V(y)}}{\int e^{-\beta V}}dy
\end{equation}
we obtain: for $t$ sufficiently large, the {\it total variation norm}
can be bounded as
\begin{equation}
  \label{e:TVconvergence}
  \norm{\mu_t - \nu}_\TV \equiv \sup_{\norm{f}\leq 1}\abs{\int f(x) \mu_t(dx)
    - \int f(x) \nu(dx)}\leq C(\mu_0) e^{-(\lambda_2 - \lambda_1) t}.
\end{equation}
In the above expression, $\norm{f}=\norm{f}_{L^\infty(\Omega)}$.  This
shows that if the process remains in $\Omega$ for a sufficiently large
amount of time (typically of the order of $1/(\lambda_2 -
\lambda_1)$), then its law at time $t$ is close to the QSD $ \nu$.

Since we are interested in ensuring that the state to state dynamics
are accurate, we observe that \eqref{e:TVconvergence} implies
agreement of the exit distribution of $(T, X_T)$ in total variation
norm between processes initially distributed according to $\mu_t$ and
$\nu$.  Indeed, starting from the probability measure $\mu_t$, given
any measurable $g:\R^+\times \partial \Omega \to \R$, we see that exit
distribution observables can be reformulated as observables on
$\Omega$ against $\mu_t$:
\begin{equation}
  \E^{\mu_t}\bracket{g(T, X_T)} = \int \underbrace{\E^{x}\bracket{g(T,
      X_T)}}_{\equiv G(x)}
  \mu_t(dx) = \int G(x) \mu_t(dx).
\end{equation}
Therefore,
\begin{equation}
  \sup_{\norm{g}\leq
    1}\abs{\E^{\mu_t}\bracket{g(T, X_T)} -\E^{\nu}\bracket{g(T, X_T)}
  }\leq \norm{\mu_t - \nu}_{\TV} \leq C(\mu_0) e^{-(\lambda_2 - \lambda_1)t},
\end{equation}
where $\norm{g} = \norm{g}_{L^\infty(\R^+\times \partial \Omega)}$.
Thus, convergence of $\mu_t $ to $\nu$ implies agreement of the exit
distributions, starting from $\mu_t$ and $\nu$.

We are now in position to introduce the original ParRep algorithm.

\subsection{Parallel Replica Dynamics}
\label{sec:orig_parrep}

The goal of the ParRep algorithm is to rapidly generate a physically
consistent first hitting point and first exit time for each visited
state.  Information about where, precisely, the trajectory is within
each state will be sacrificed to more rapidly obtain this information.

In the following, we assume that we are given a partition of the
configuration space $\R^d$ into states, and we denote by $\Omega$ one
generic element of this partition. We also assume that we have $N$
CPUs available for parallel computation.

The original ParRep algorithm~\cite{Voter:1998p13729} is implemented
in three steps, repeated as the process moves from one state to
another. It requires the specification, {\it a priori}, of two times
to equilibrate to each state, $\tcorr$ and $\tphase$.  Let us consider
a single reference process, $X_t^\refe$, with $X_0^\refe \sim \mu_0$
evolving under \eqref{e:odlang}, and set the simulation clock,
corresponding to the physical time, to zero, $t_{\simu}=0$.  The
simulation clock $t_{\simu}$ will be updated during the algorithm.
\begin{description}

\item[Decorrelation Step] Let $\Omega$ denote the state in which
  $X_{t_{\simu}}^\refe$ currently resides. If the trajectory has not
  left $\Omega$ after running for $\tcorr$ amount of time, the
  algorithm proceeds to the dephasing step, the simulation clock being
  advanced as
  \begin{equation*}
    t_{\simu}\mapsto t_{\simu} + \tcorr.
  \end{equation*}
  Otherwise a new decorrelation starts from the new state, the
  simulation clock being advanced as
  \begin{equation*}
    t_{\simu}\mapsto t_{\simu} +  T^\refe,\quad T^\refe =
    \inf\set{t\geq 0\mid X^\refe_{t_{\simu}+t} \notin \Omega}.
  \end{equation*}

\item[Dephasing Step] In this step, $N$ independent and identically
  distributed samples of the QSD of $\Omega$ are generated.  These
  samples will be distributed over the $N$ CPUs and will be used as
  initial conditions in the subsequent parallel step. During the
  dephasing step, the counter $t_{\simu}$ is not advanced.

\noindent
In the original ParRep, the sampling of the QSD is accomplished by a
rejection algorithm. For $k=1,\ldots, N$, generate a starting point
$\tilde X_0^k\sim \eta_0$ and integrate it under~\eqref{e:odlang}
until either time $\tphase$ or $\tilde X_t^k$ leaves $\Omega$. Here,
$\eta_0$ denotes any distribution with support in $\Omega$ (for
example, a Dirac mass at the end point of the reference trajectory,
after the decorrelation step).  If $\tilde X_t^k$ has not exited
before time $\tphase$, set the $k$-th replica's starting point $X_0^k
= \tilde X_{\tphase}^k$ and advance $k\mapsto k+1$. Otherwise, reject
the sample, and start a new trajectory with $\tilde X_0^k \sim
\eta_0$.  Since these samples are independent, they can be generated
in parallel.

\item[Parallel Step] Let the $N$ samples obtained after the dephasing
  step evolve under \eqref{e:odlang} in parallel (one on each CPU),
  driven by independent Brownian motions, until one escapes from
  $\Omega$. Let us denote
$$k_\star = {\rm
  argmin}_k T^k$$ the index of the first replica which exits $\Omega$.
During the time interval $[t_{\simu},t_{\simu} + N T^{k_\star}]$, the
reference process is defined as trapped in $\Omega$. Accordingly, the
simulation clock is advanced as
\begin{equation}
  \label{e:tadvance_par}
  t_{\simu}\mapsto t_{\simu} + N T^{k_\star}, \quad 
\end{equation}
The first replica to escape becomes the new reference process. A new
decorrelation step now starts, applied to the new reference process
with starting point $X^\refe_{t_{\simu}}=X^{k_\star}_{T^{{k_\star}}}$.
\end{description}

The justifications underlying the ParRep algorithm are the
following. Using the third property~\eqref{e:TVconvergence} of the
QSD, it is clear that if $\tcorr$ is chosen sufficiently large, then,
at the end of the decorrelation step, the reference process is such
that $X^\refe_{t_{\simu}}$ is approximately distributed according to
the QSD. This same property explains why the rejection algorithm used
in the dephasing step yields (approximately) $N$ i.i.d. samples
distributed according to the QSD, at least if $\tphase$ is chosen
sufficiently large. Finally the second property~\eqref{e:qsd_indep}
justifies the parallel step; since the $N$ replicas are i.i.d. and
drawn from the QSD, they have exponentially distributed exit times and
thus $\prob^\nu\bracket{\min\set{T^1,\ldots, T^N}>t} = \prob^\nu[T^1>
Nt]$. Moreover, by the independence property in~\eqref{e:qsd_indep},
the exit points $X^{k_\star}_{T^{{k_\star}}}$ and $X^{1}_{T^{1}}$ have
the same distribution.

Notice that it is the magnification of the first exit time by a factor
of $N$ in the parallel step that yields the speedup in terms of wall
clock time. If the partition of the configuration space is chosen in
such a way that, most of the time, the stochastic process exits from
the state before having reached the QSD (namely before $\tcorr$),
there is no speedup. In this case, ParRep essentially consists in
following the reference process. There is no error, but no gain in
performance, and computational resources are wasted. To observe a
significant speedup, the partition of the configuration space should
be such that most of the defined states are metastable, in the sense
that the typical exit time from the state is much larger than the time
required to approximate the QSD.

Of course, the QSD is only sampled approximately, and this introduces
error in ParRep.  The time $\tcorr$ (resp. $\tphase$) must be
sufficiently large such that $\prob^{\mu_0}\bracket{X_{\tcorr}^{\refe}
  \in \bullet\mid T^{\refe} >\tcorr}\approx \nu$ (resp. such that
$\prob^{\eta_0}\bracket{X_{\tphase} \in \bullet\mid T >\tphase}\approx
\nu$).  The mismatch between the distributions at times $\tcorr$ and
$\tphase$, directly, and independently, contribute to the overall
error of ParRep; see \cite{Simpson:2013cs}. Also note that these
parameters are state dependent.  In view of~\eqref{e:TVconvergence},
one may think that a good way to choose $\tcorr$ and $\tphase$ is to
consider a multiple of $1/(\lambda_2 - \lambda_1)$. This is
unsatisfactory for two reasons. First, it is difficult to numerically
compute the spectral gap $\lambda_2 - \lambda_1$ because of the
high-dimensionality of the associated elliptic problem.  Second, the
pre-factors $C(\mu_0)$ and $C(\eta_0)$ in \eqref{e:TVconvergence} are
also difficult to evaluate, and could be large.

In view of the preceding discussion, ParRep can be applied to a wide
variety of problems and for any predefined partition of the
configuration space into states provided one has:
\begin{itemize}
\item A way to construct an adequate $\tcorr$ (or more precisely to
  assess the convergence of $\mu_t$ to $\nu$) for each state;
\item A way to sample the QSD of each state.
\end{itemize}
The aim of the next section is to provide a modified ParRep algorithm
to deal with these two difficulties.

\section{The Modified ParRep Algorithm}
\label{sec:modParRep}

We propose to use a branching and interacting particle system (the
Fleming-Viot particle process) together with convergence diagnostics
to simultaneously and dynamically determine adequate values $\tcorr$
and $\tphase$, while also generating an ensemble of $N$ samples from a
distribution close to that of the QSD.

\subsection{The Fleming-Viot Particle Process}
\label{s:fvprocess}

In this section, we introduce a branching and interacting particle
system which will be one of the ingredients of the modified ParRep
algorithm. This process is sometimes called the Fleming-Viot particle
process~\cite{ferrari07:_quasi_flemin_viot}.

\begin{figure}
  \subfigure[$t_1$]{\includegraphics[width=4cm]{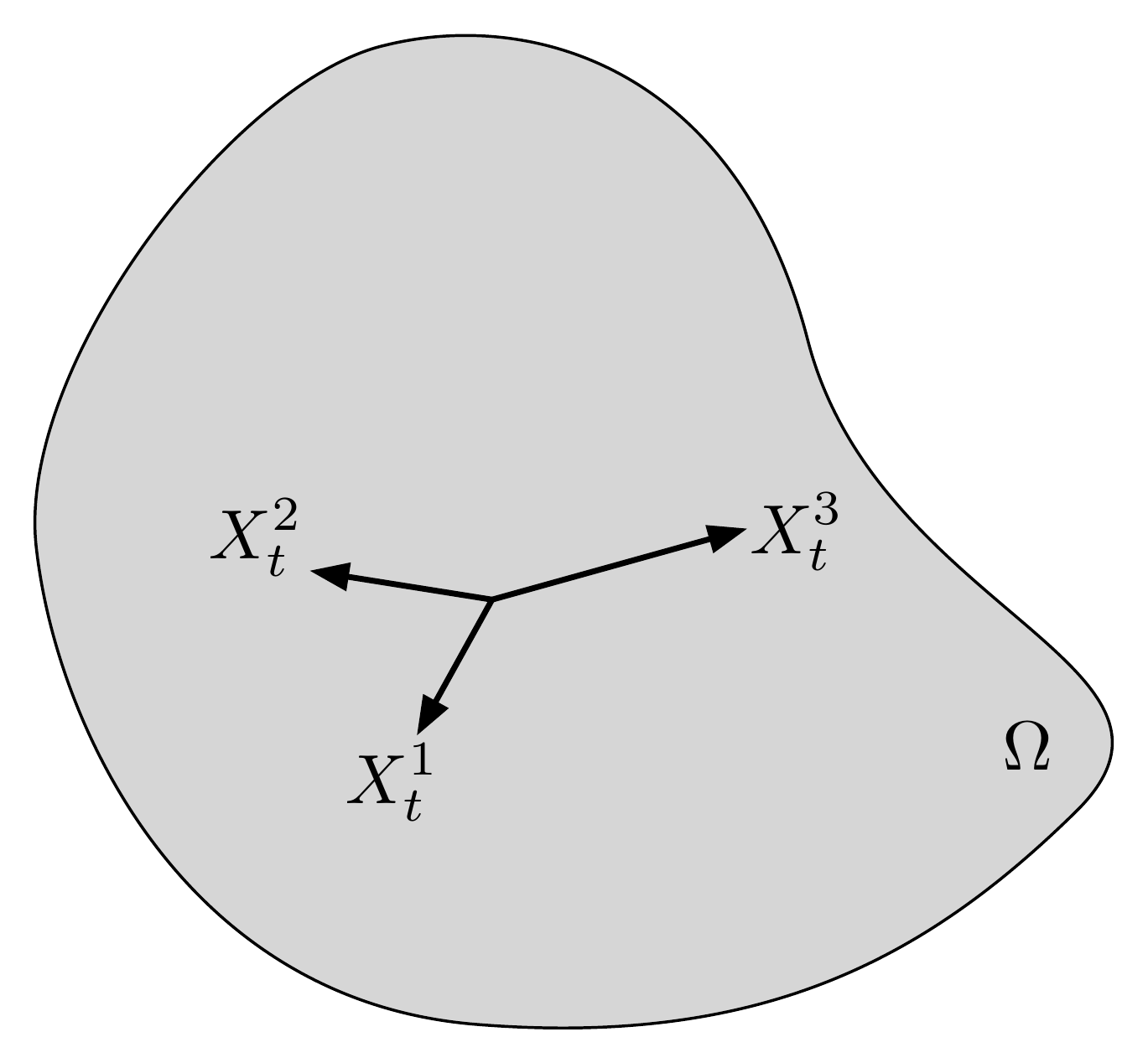}}
  \subfigure[$t_2$]{\includegraphics[width=4cm]{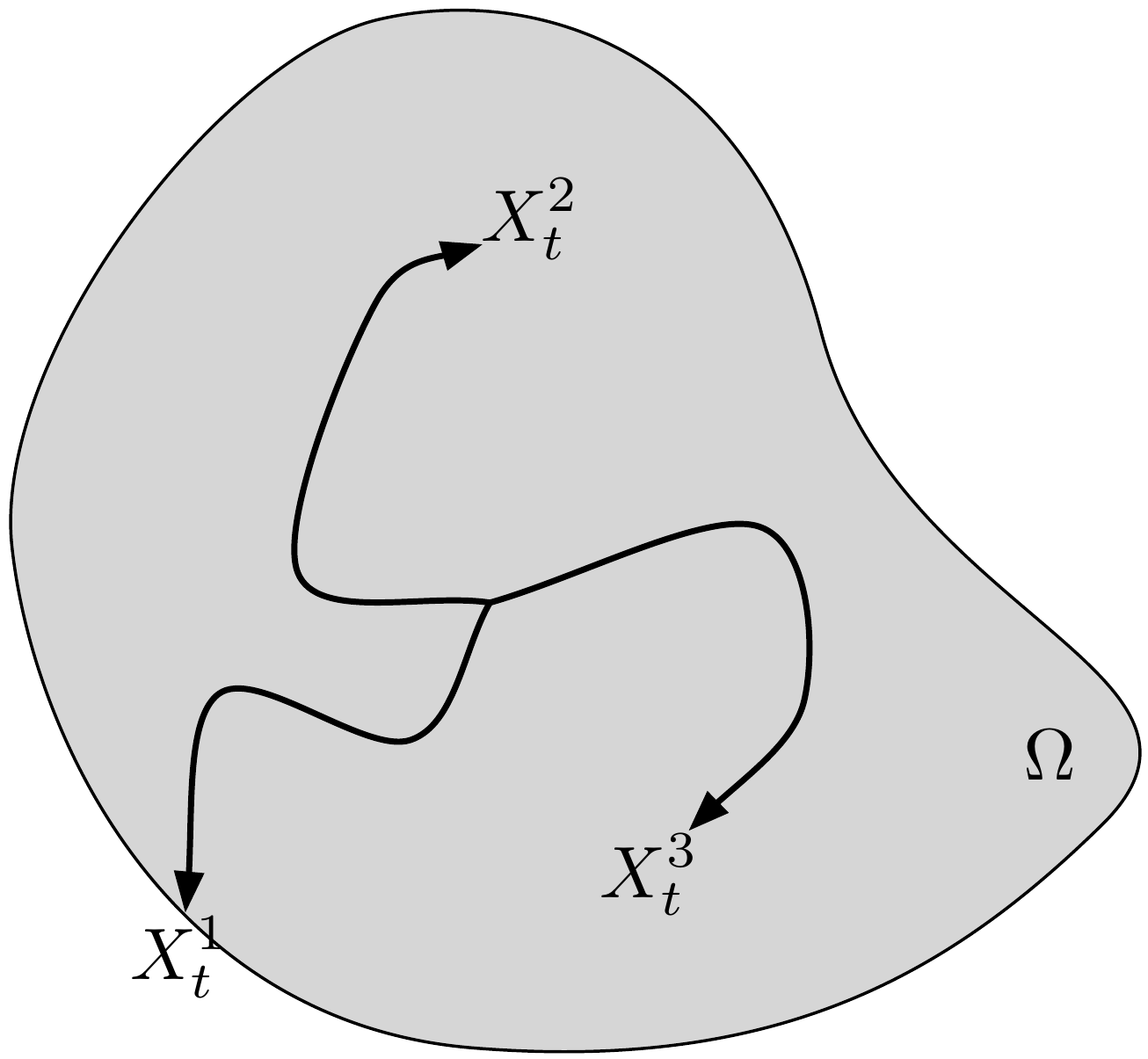}}

  \subfigure[$t_2$]{\includegraphics[width=4cm]{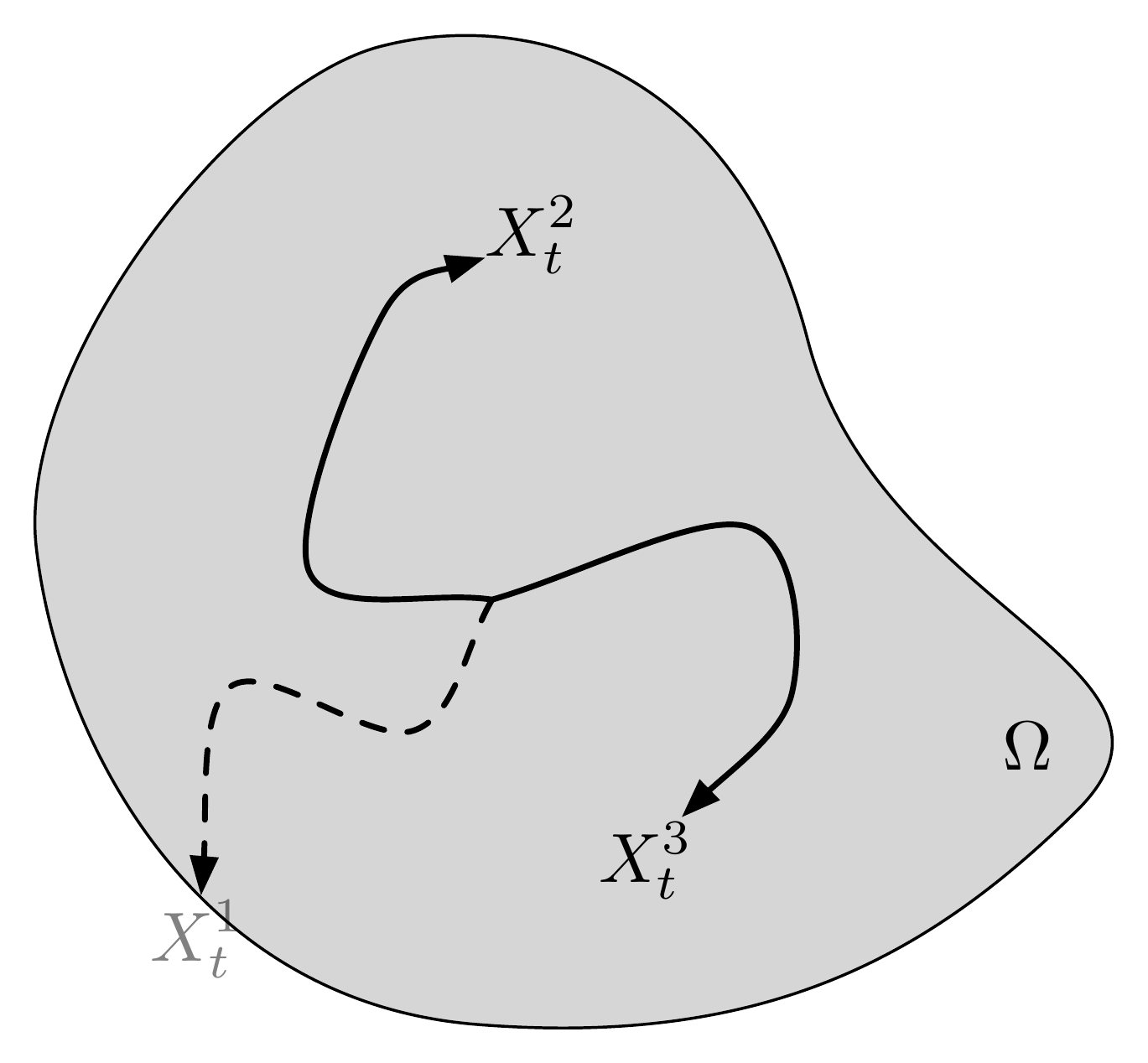}}
  \subfigure[$t_3$]{\includegraphics[width=4cm]{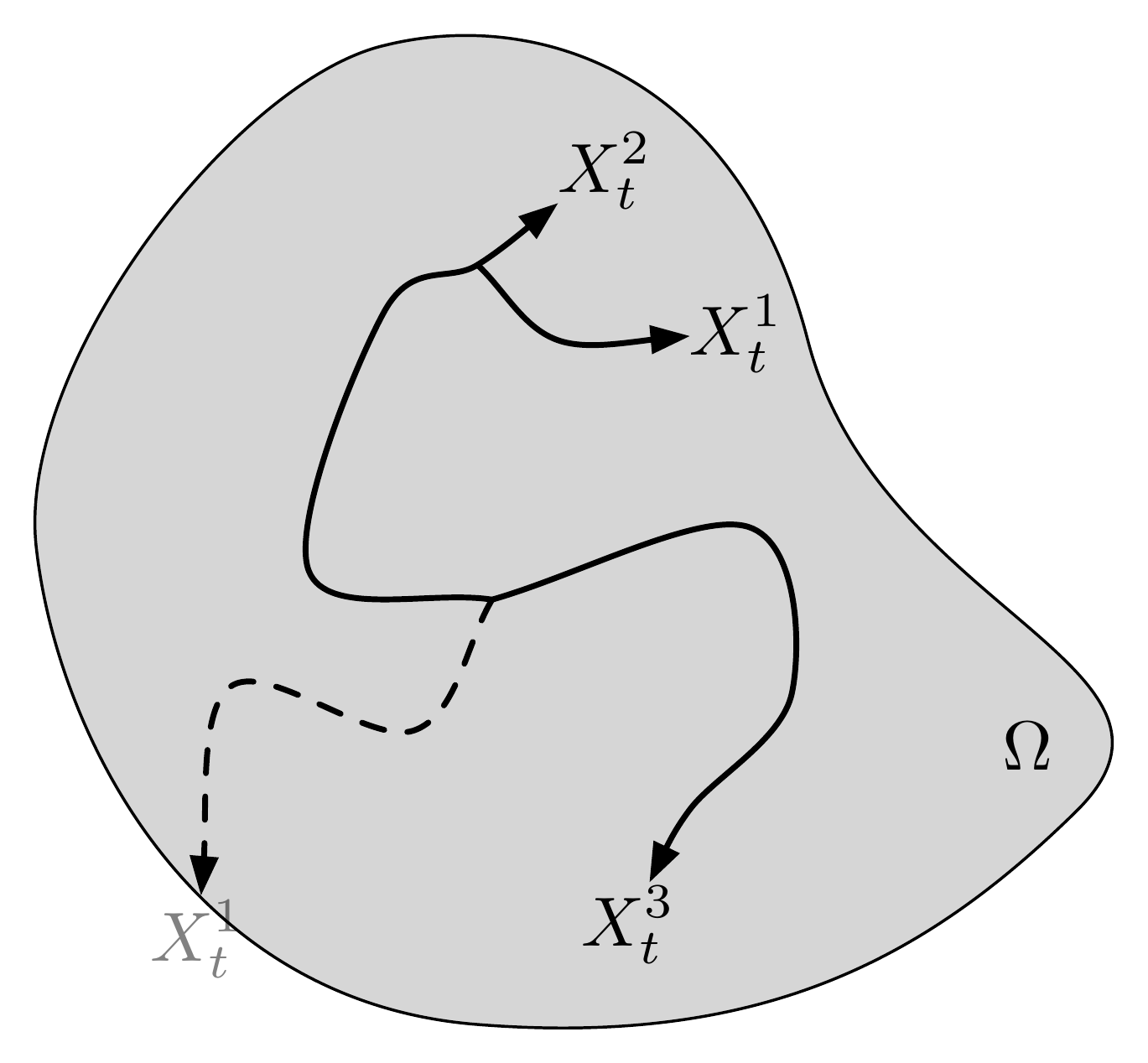}}
  \caption{The branching \& interacting particle system used to sample
    the QSD in the case $N=3$ at three times: $t_1<t_2<t_3$.  The
    trajectories run independently until one exits, as in (b).  The
    process that has reached the boundary is killed, as in (c).  Then
    a survivor is instantaneously branched to maintain a constant
    number of trajectories, as in (d).}
  \label{f:bips_diagram}
\end{figure}

Let us specify the Fleming-Viot particle process; see also the
illustration in Figure \ref{f:bips_diagram}. Let us consider
i.i.d. initial conditions $X_0^k$ ($k \in \{1, \ldots,N\}$)
distributed according to $\mu_0$, a probability distribution with
support in $\Omega$. The process is as follows:
\begin{enumerate}
\item Integrate $N$ realizations of \eqref{e:odlang} with independent
  Brownian motions until one of them, say $X_t^1$, exits;
\item Kill the process that exits;
\item With uniform probability $1/(N-1)$, randomly choose one of the
  survivors, $X_t^2,\ldots, X_t^N$, say $X_t^{2}$;
\item Branch $X_t^{2}$, with one copy persisting as $X_t^{2}$, and the
  other becoming the new $X_t^1$ (and thus evolving in the future
  independently from $X_t^2$).
\end{enumerate}
We denote this branching and interacting particle process by ${\bf
  X}_{t,N}=(X_t^1, \ldots,X_t^N)$, and define the associated empirical
distribution
\begin{equation}
  \mu_{t,N} \equiv \frac{1}{N} \sum_{k=1}^N \delta_{X^k_t}.
\end{equation}
The Fleming-Viot particle process can be implemented in parallel, with
each replica $X_t^k$ evolving on distinct CPUs. The communication cost
(due to the branching step) will be small, provided the state under
consideration is such that the exit events are relatively rare; {\it
  i.e.}, it is metastable.

The Fleming-Viot particle process has been studied for a variety of
underlying stochastic processes; see, for example,
\cite{ferrari07:_quasi_flemin_viot,Meleard:2012vl,del2004feynman} and
the references therein.  In \cite{Meleard:2012vl}, the authors prove
that for a problem in dimension one, the following relation holds: for
any $A \subset \Omega$,
\begin{equation}
  \label{e:fvconv1}
  \lim_{N\to \infty} \mu_{t,N}(A)  = 
  \mu_t(A).
\end{equation}
From~\eqref{e:TVconvergence} and~\eqref{e:fvconv1}, we infer that
$\lim_{t\to \infty} \lim_{N\to \infty} \mu_{t,N}(A) = \nu (A)$.  This
result is anticipated to hold for general dynamics, including \eqref
{e:odlang}.

The property~\eqref{e:fvconv1} of the Fleming-Viot particle process is
instrumental in our modified ParRep algorithm. It will be used in two
ways:
\begin{itemize}
\item Since an ensemble of realizations distributed according to
  $\mu_t$ can be generated using the Fleming-Viot particle process, we
  will assess convergence of $\mu_t$ to the stationary distribution
  $\nu$ by applying convergence diagnostics to the ensemble $(X^1_t,
  \ldots ,X^N_t)$. This will give a practical way to estimate the time
  required for convergence to the QSD in the decorrelation step, by
  {\it simultaneously} running the decorrelation step (on the
  reference process) and the dephasing step using a Fleming-Viot
  particle process on $N$ other samples, starting from the same
  initial condition as the reference process: the decorrelation time
  is estimated as the time required for convergence to a stationary
  state for the Fleming-Viot particle process.
\item In addition, the Fleming-Viot particle process introduced in the
  procedure described above gives a simple way to sample the QSD.  We
  use the replicas generated by the modified dephasing step at the
  time of stationarity.
\end{itemize}

The modified ParRep algorithm will be based on the two following
assumptions on the Fleming-Viot particle process.  While we do not
make this rigorous, we believe it could be treated in specific cases,
and our numerical experiments show consistency with direct numerical
simulation.
\begin{description}
\item[Assumption (A1)] For sufficiently large $N$, $\mu_{t,N}$ is a
  good approximation of $\mu_t$;
\item[Assumption (A2)] The realizations generated by the Fleming-Viot
  particle process are sufficiently weakly correlated so as to allow
  the use of both the convergence diagnostics presented below and the
  temporal acceleration expression~\eqref{e:tadvance_par}, which both
  assume independence.
\end{description}

As already mentioned above, the first assumption is likely satisfied
in our setting, though we were not able to find precisely this result
in the literature.  See~\cite{Meleard:2012vl} for such a result in a
related problem.

The second assumption is more questionable. We make two comments on
this. First, our numerical experiments show that the modified ParRep
algorithm (which is partly based on (A2)) indeed yields correct
results compared to direct numerical simulation; thus, the assumption
is not grossly wrong, at least in these settings. Second, the
correlations introduced by the Fleming-Viot particle process are most
likely a concern for problems where the state is only weakly
metastable. In truly metastable states, exits will be infrequent, so
the correlations amongst the replicas will be weak.  For states which
are not metastable, the reference process will likely exit before
stationarity can be reached, rendering the concern moot.  It is
therefore in problems between the two cases that practitioners may
have some cause for concern.

There are several ways to ameliorate reservations about the second
assumption. First, it is known that for such branching and interacting
particle systems, a propagation of chaos result holds,
\cite{sznitman-91}.  This means that if we run Fleming-Viot with $M\gg
N$ processes, then, as $M$ tends to infinity, with $N$ fixed, the
first $N$ trajectories in the process become i.i.d.  Second, one could
run a separate Fleming-Viot particle process for each of the $N$
replicas, retaining only the first trajectory of each of the $N$
i.i.d. Fleming-Viot particle processes.  Finally, the Fleming-Viot
particle process, together with convergence diagnostics, could be used
to identify appropriate values of $\tphase$ and and $\tcorr$, as
explained below in Section \ref{s:diags}. This value $\tphase$ could
then be used in the original rejection sampling algorithm, run in
tandem with the Fleming-Viot particle process.  This would provide
independent samples from the QSD.

\subsection{Convergence Diagnostics}
\label{s:diags}

While the Fleming-Viot particle process gives us a process that will
converge to the QSD, there is still the question of how long it must
be run in order for $\mu_{t,N}$ (and thus $\mu_t$ according to (A1))
to be close to equilibrium.  This is a ubiquitous problem in applied
probability and stochastic simulation: when sampling a distribution
via Markov Chain Monte Carlo, how many iterations are sufficient to be
close to the stationary distribution?  For a discussion on this
general issue, see for example,
\cite{Brooks:1998ux,brooks2011handbook,Cowles:1996vt}.  We propose to
use convergence diagnostics to test for the stationarity of
$\mu_{t,N}$.  When a user specified convergence criterion is
satisfied, $\mu_{t,N}$ is declared to be at its stationary value and
the time at which this occurs is taken to be $\tcorr$ (and $\tphase$).

We have found Gelman-Rubin statistics to be effective for this
purpose, \cite{Gelman:1992ts,Brooks:1998um,Brooks:1998ux}.  In the
simplest form, such statistics compute the ratio of two estimates of
the asymptotic variance of a given observable.  Since numerator and
denominator estimate the same quantity, the ratio converges to one.

The statistic can be defined as follows.  Let $\obs:\Omega\to \R$ be
some observable, and let
\begin{equation}
  \bar{\obs}^k_{t} \equiv t^{-1} \int_0^t \obs(X_s^k)ds, \quad \bar{\obs}_{t}\equiv \frac{1}{N} \sum_{k=1}^N \bar{\obs}^k_{t}
  = \frac{1}{N} \sum_{k=1}^N t^{-1} \int_0^t \obs(X_s^k)ds,
\end{equation}
be the average of an observable along each trajectory and the average
of the observable along all trajectories.  Then the statistic of
interest for observable $\obs$ is
\begin{equation}
  \label{e:R2}
  \hat{R}_{t}(\obs) = \frac{\frac{1}{N}\sum_{k=1}^N t^{-1} \int_0^t(\obs(X_s^k) - \bar{\obs}_{t})^2 ds }{\frac{1}{N}
    \sum_{k=1}^N t^{-1} \int_0^t (\obs(X_s^k) -\bar{\obs}_{t}^k )^2 ds }.
\end{equation}
Notice that $\hat R_t(\obs) \geq 1$, and as all the trajectories
explore $\Omega$, $\hat{R}_t(\obs)$ converges to one as $t$ goes to
infinity.  {We also observe that since \eqref{e:R2} is a
  trajectory average, both its bias and variance will be
  $\bigo(t^{-1})$, \cite{Asmussen:2010wy}.}

These statistics were not developed with the intention of handling
branching interacting particle systems.  The authors had in mind that
the $N$ trajectories would be completely independent, which is not the
case for the Fleming-Viot particle process. This is one reason why we
introduced Assumption (A2) above.  However, we will demonstrate in the
numerical experiments below that this convergence diagnostic indeed
provides meaningful results for the Fleming-Viot particle process (see
in particular Section~\ref{sec:per1d}).

Here, we caution the reader that all convergence diagnostics are
susceptible to the phenomena of {\it pseudo-convergence}, which occurs
when a particular observable or statistic appears to have reached a
limiting value, and yet the empirical distribution of interest remains
far from stationarity; see \cite{brooks2011handbook}.  This can occur,
for instance, if the state has an internal barrier that obstructs the
process from migrating from one mode to the other. A computational
example of this is given below, in Section~\ref{sec:DP2D}.

There is still the question of what observables to use in computing
the statistics.  Candidates include:
\begin{itemize}
\item Moments of the coordinates;
\item Energy $V(x)$;
\item Distances to reference points in configuration space.
\end{itemize}
Assuming they are not costly to evaluate, as many such observables
should be used; see the example in Section~\ref{sec:DP2D}.

Our test for stationarity is as follows.  Given some collection of
observables $\{\obs_j:\R^d \to \R\}_{ j \in \{1, \ldots, J\}}$, their
associated statistics $\{\hat{R}_t(\obs_j)\}_{ j \in \{1, \ldots,
  J\}}$, and a tolerance $\tol >0$, we take as a stationarity
criterion:
\begin{equation}
  \label{e:termination}
  \forall j \in \{1,
  \ldots, J\}, \, \hat{R}_t(\obs_j) < 1 +  \tol.
\end{equation}
In other words, the dephasing and decorrelation times are set as
\begin{equation}
  \label{e:termination_times}
  \tphase =\tcorr= \inf \set{t\geq 0\mid \hat{R}_t(\obs_j) < 1 +  \tol, \;
    \forall j}.
\end{equation}

\subsection{The Modified ParRep Algorithm}

We now have the ingredients needed to present the modified ParRep
algorithm (which should be compared to the original ParRep given in
Section~\ref{sec:orig_parrep}).  Let us consider a single reference
process, $X_t^\refe$, with $X_0^\refe \sim \mu_0$ evolving under
\eqref{e:odlang}, and let us set $t_{\simu}=0$.

\begin{description}

\item[Decorrelation and Dephasing Step] Denote by $\Omega$ the state
  in which $X_{t_{\simu}}^\refe$ lives. The decorrelation and
  dephasing steps are carried out at the same time, in parallel: the
  reference process $X_t^\refe$ and the Fleming-Viot particle process
  ${\bf X}_{t,N}$ begin at the same time from the same point in
  $\Omega$.  Convergence diagnostics are assessed on $\mu_{t,N}$, and
  when the stationarity criterion~\eqref{e:termination} is satisfied,
  in the case that the reference process has never left $\Omega$, both
  decorrelation and dephasing steps are terminated, and one proceeds
  to the Parallel Step, after advancing the simulation clock as
  \begin{equation*}
    t_{\simu}\mapsto t_{\simu} + \tcorr,
  \end{equation*}
  $\tcorr$ being defined by~\eqref{e:termination_times}.  In this
  case, the decorrelation/dephasing step is said to be successful.

  If at any time before reaching stationarity the reference process
  leaves $\Omega$, the Fleming-Viot particle process terminates, ${\bf
    X}_{t,N}$ is discarded, the simulation clock is advanced as
  \begin{equation*}
    t_{\simu}\mapsto t_{\simu} +  T^\refe,\quad T^\refe =
    \inf\set{t\geq 0\mid X^\refe_{t_{\simu}+t} \notin \Omega}.
  \end{equation*}
  The process $X_t^{\refe}$ then proceeds into the new state, where a
  new decorrelation/dephasing step starts. In this case, the
  decorrelation/dephasing step is said to be unsuccessful.

\item[Parallel Step] The parallel step is similar to the original
  parallel step. Consider the $N$ positions of ${\bf X}_{\tcorr,N}$
  obtained at the end of the dephasing step as initial
  conditions. These are then evolved in parallel
  following~\eqref{e:odlang}, driven by independent Brownian motions,
  until one replica, with index $k_\star$, escapes from $\Omega$. The
  simulation clock is then advanced according to
  \eqref{e:tadvance_par},
  \[
  t_{\simu}\mapsto t_{\simu} + N T^{k_\star}
  \]
  The replica which first exits becomes the new reference process, and
  a new decorrelation/dephasing step starts.
\end{description}

{Our modified ParRep algorithm differs from the original
  in two essential ways, which merit comment. First, as already
  mentioned, the main attraction of the modified ParRep algorithm is
  that the convergence times $\tcorr$ and $\tphase$ do not need to be
  chosen {\it a priori}, but are instead computed on the fly using
  physically informed observables. This is why the Fleming-Viot
  particle process (combined with a convergence diagnostic) is
  essential to our algorithm. This is the major improvement of
  modified ParRep, which broadens the applicability of the algorithm
  to a general partition of the configuration space, as {\it a priori}
  estimates for $\tcorr$ and $\tphase$ are, in general, unlikely to be
  available.

  Second, in the modified ParRep, the Fleming-Viot particle process is
  also used to sample the QSD (the dephasing step). Part of the appeal
  of the Fleming-Viot particle process is that it is robust enough to
  sample states which are not strongly metastable (the typical exit
  time is not dramatically larger than the time to converge to the
  QSD).  For highly metastable regions, rejection sampling and the
  Fleming-Viot particle process will yield similar results, but in
  more general scenarios, this will not be the case.  Indeed, while
  the Fleming-Viot particle process is well defined in the limit of
  time tending to infinity, rejection sampling will degenerate since
  all replicas will eventually exit. This is critical when applying
  convergence diagnostics, as the time to stop is found dynamically by
  condition \eqref{e:termination}.  Therefore, if one requires $N$
  replicas distributed according to the QSD before proceeding to the
  parallel step, the Fleming-Viot particle process may, therefore, be
  more efficient than rejection sampling.

  Our discussion of the Fleming-Viot particle process, as a QSD
  sampling strategy, would not be complete without two
  comments. First, some implementations of the rejection algorithm in
  ParRep do not require all $N$ replicas to be dephased before
  proceeding to the parallel step.  Replicas are run asynchronously,
  and as soon as one has reached the {\it a priori} value of
  $\tphase$, it is immediately promoted to the parallel step.  This is
  how the algorithm is implemented in, for instance,
  \cite{Kum:2004iu}, and it should be taken into account when
  comparing QSD sampling algorithms.  Second, as already mentioned
  above, the robustness of the Fleming-Viot process comes at the cost
  of generating correlated samples, which are not easy to control.
  Two methods for overcoming correlations amongst the replicas are
  proposed at the end of Section \ref{s:fvprocess}.}

We will illustrate this modified ParRep algorithm on various numerical
examples in Section~\ref{s:parrep_examples}, after a discussion of the
performance and limitations of the convergence diagnostics in
Section~\ref{s:offline_examples}.

\section{Illustration of Convergence Diagnostics}
\label{s:offline_examples}
In this section we present two numerical examples showing the
subtleties of the Gelman-Rubin statistics and the broader problems
raised by stationarity testing.  These are ``offline'' in the sense
that they are not used as part of the ParRep algorithm here.  They
show that the Gelman-Rubin statistics are consistent with our
expectations, but also susceptible to pseudo-convergence.

In both of these examples, $N=10^4$ replicas are used, and the
stochastic differential equation~\eqref{e:odlang} is discretized using
Euler-Maruyama with a time step $\Delta t = 10^{-4}$.  The Mersenne
Twister algorithm is used as a pseudo-random number generator in these
two examples, as implemented in~\cite{Galassi:2013tz}.

\subsection{Periodic Potential in 1D}
\label{sec:per1d}
For the first example, consider the 1D periodic potential $V(x) = -2
\cos(\pi x)$ at $\beta =1$ and the state $\Omega = (-1,1)$. The
initial condition is $X_0 = 0.99$.  Running the Fleming-Viot particle
process algorithm, we examine the Gelman-Rubin statistics for the
observables:
\begin{equation}
  \label{e:per1d_obs}
  x, \quad V(x), \quad \abs{x- x_{\refe}}
\end{equation}
where $x_{\refe}$ is the local minima of the current basin; $x_{\refe}
=0$ in this case.

\begin{figure}
  \includegraphics[width=8cm]{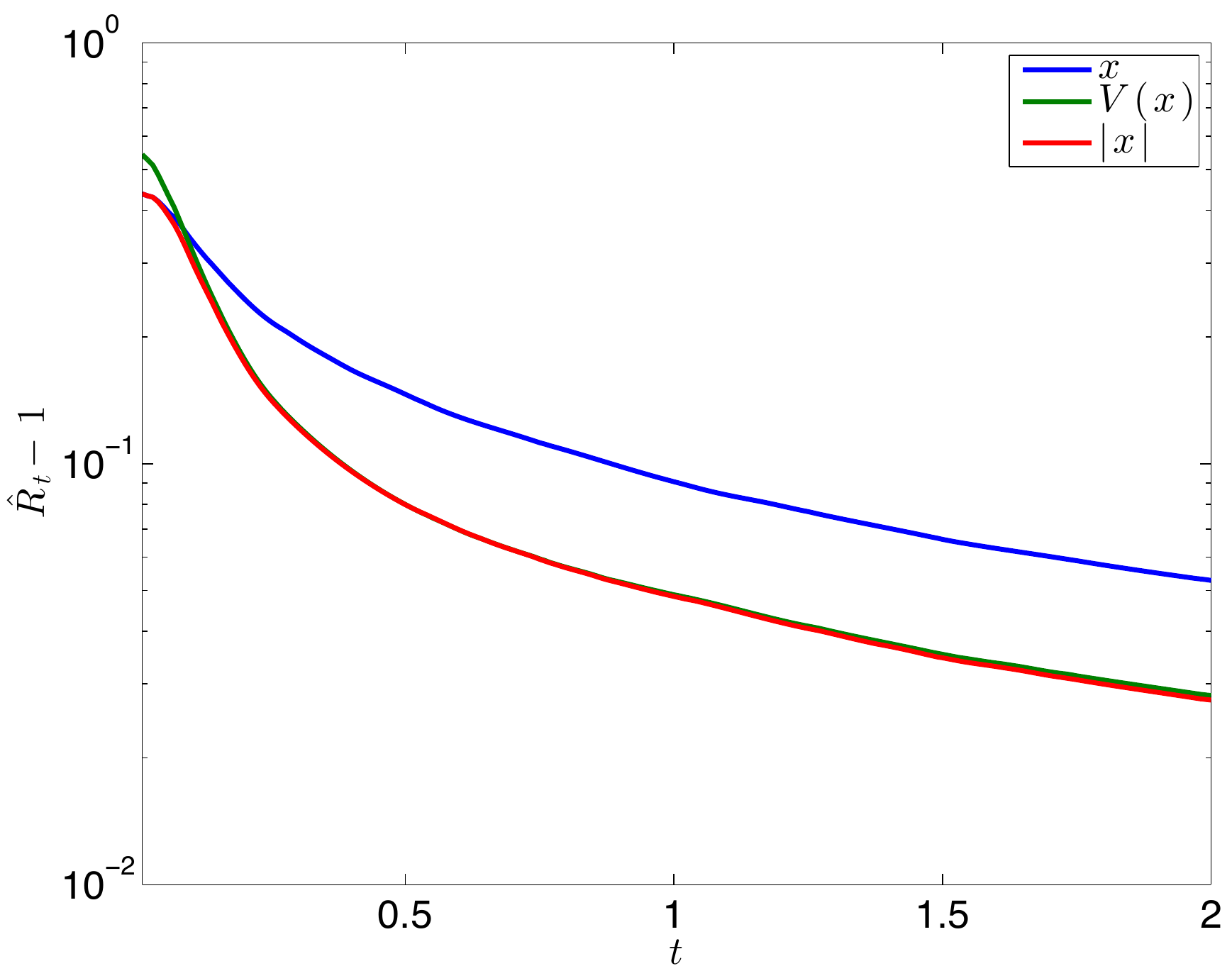}
  \caption{The Gelman-Rubin statistics as a function of time for the
    observables~\eqref{e:per1d_obs} for the Fleming-Viot particle
    process. The potential is $V(x)=-2 \cos(\pi x)$, the state is
    $(-1,1)$ and the number of replicas is $N=10^4$.}
  \label{f:per1d_grobs}
\end{figure}

\begin{figure}
  \subfigure{\includegraphics[width=6cm]{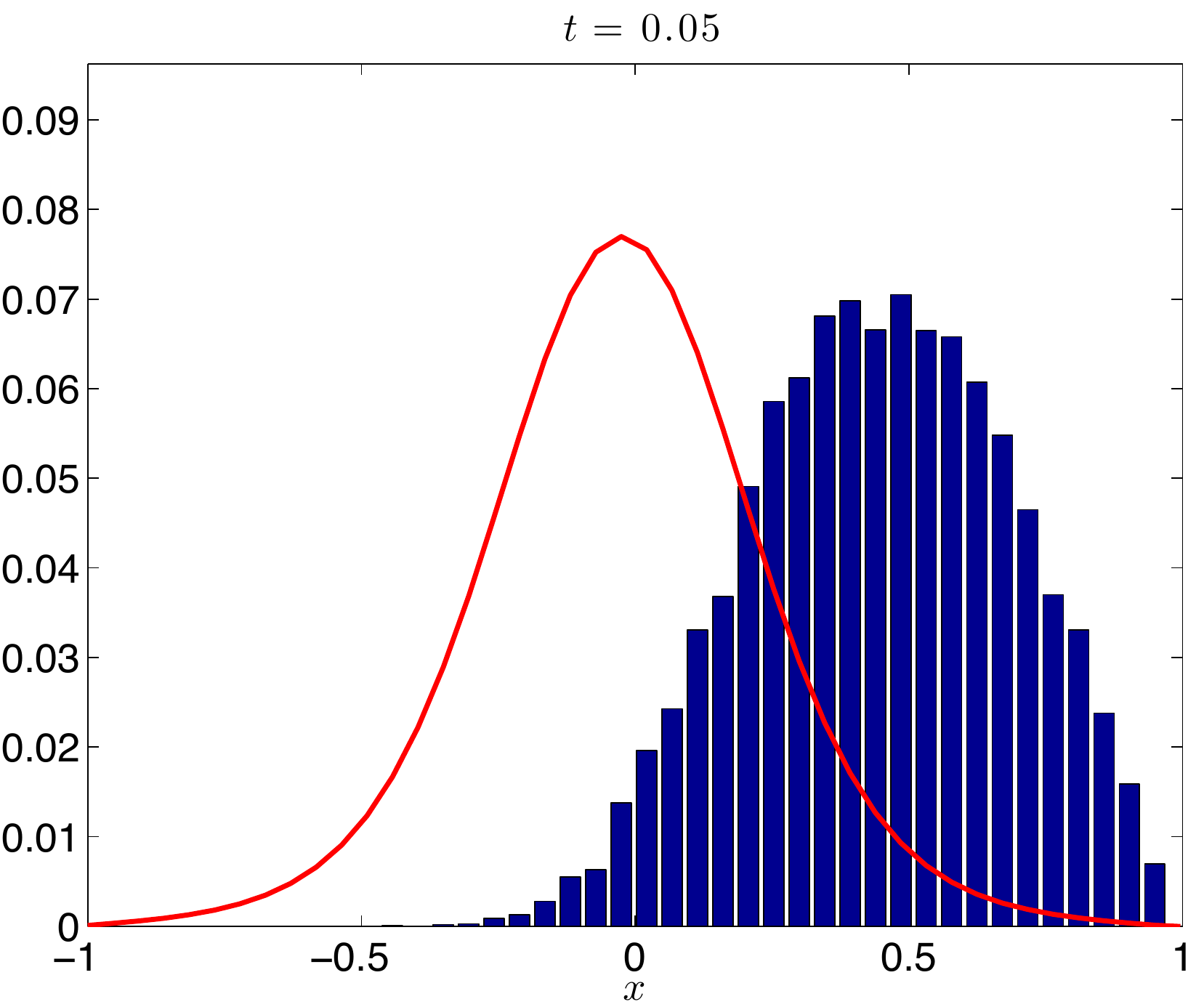}}
  \subfigure{\includegraphics[width=6cm]{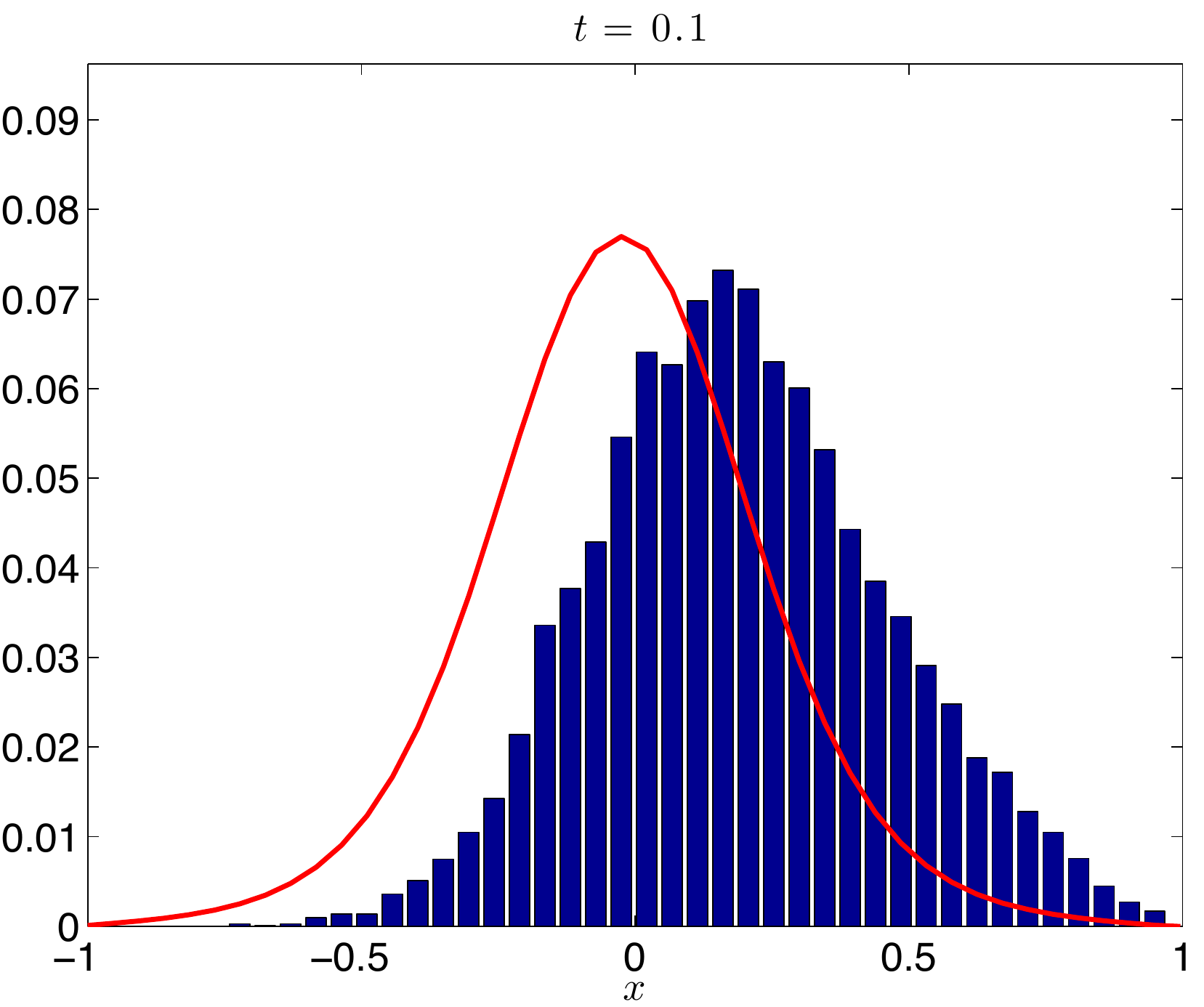}}

  \subfigure{\includegraphics[width=6cm]{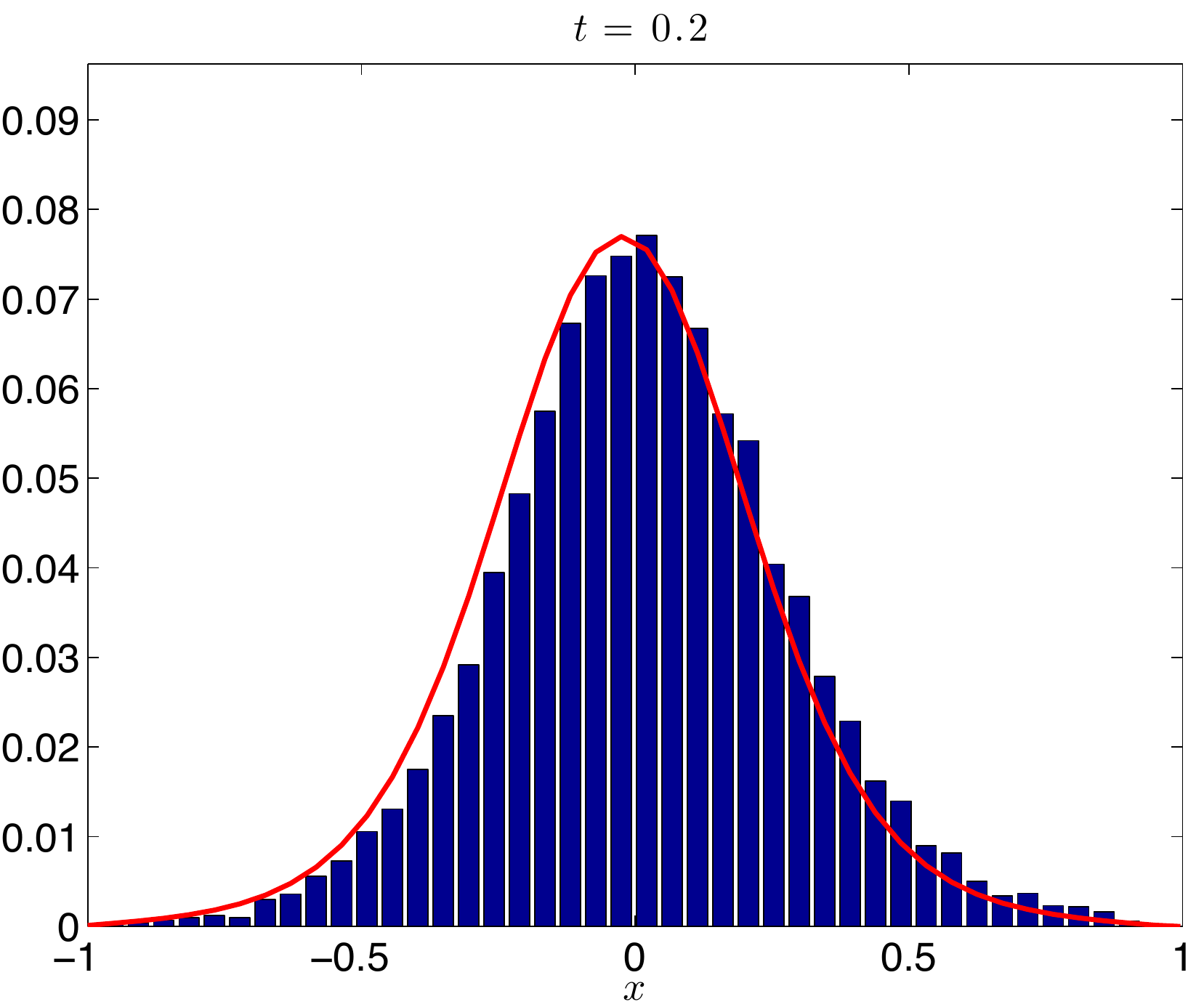}}
  \subfigure{\includegraphics[width=6cm]{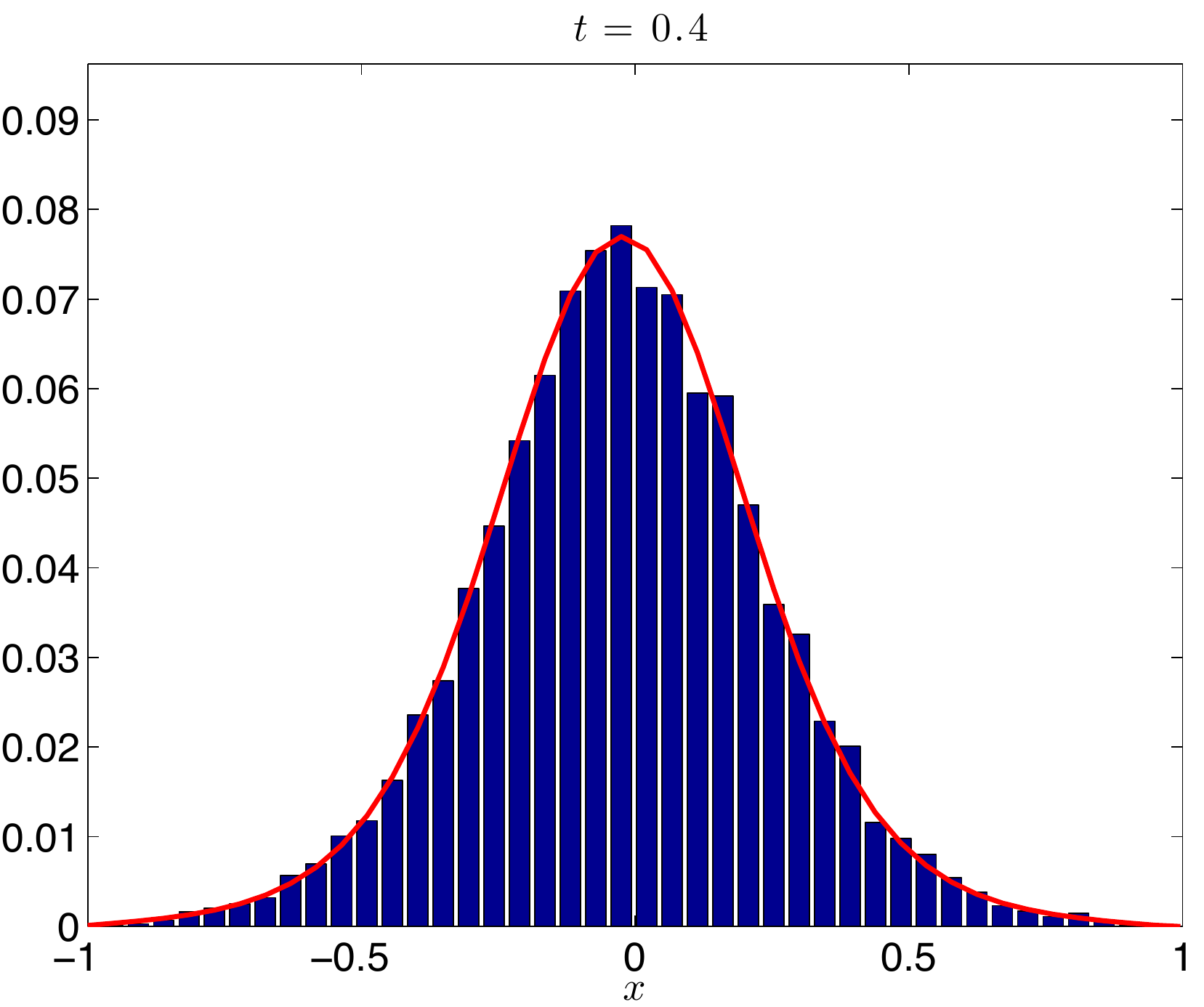}}

  \subfigure{\includegraphics[width=6cm]{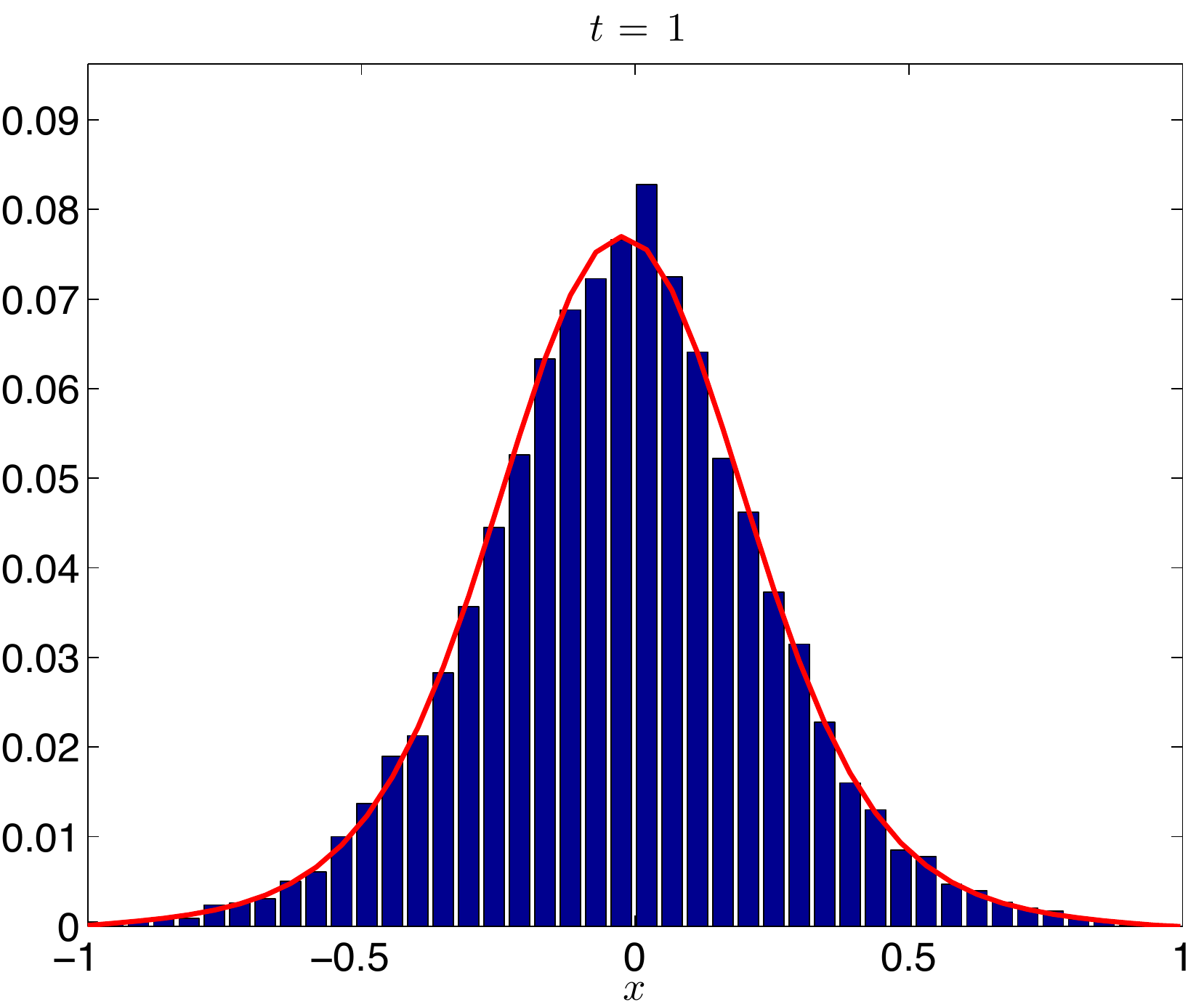}}
  \subfigure{\includegraphics[width=6cm]{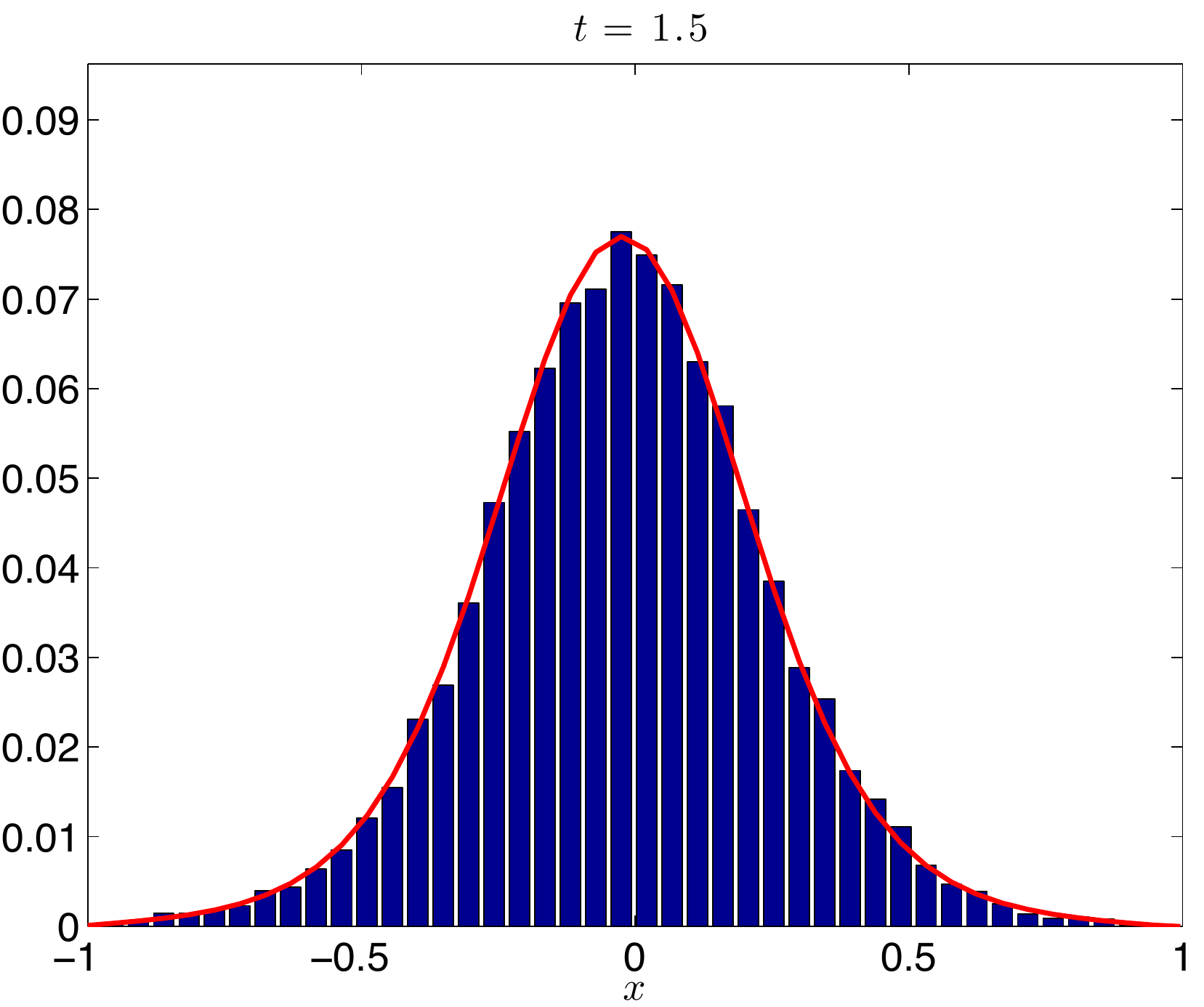}}

  \caption{Empirical distributions for the Fleming-Viot particle
    process on the 1D periodic potential $V(x)=-2 \cos(\pi x)$ over
    the domain $(-1,1)$ with $N=10^4$ replicas, see also~Figure
    \ref{f:per1d_grobs}. The red curve is the reference density of the
    QSD.}
  \label{f:per1d_hist}
\end{figure}

The Gelman-Rubin statistics, as a function of time, appear for several
observables in Figure \ref{f:per1d_grobs}.  As expected, the
statistics tend to one as time goes to infinity.  They also display
the aforementioned $\bigo(t^{-1})$ bias.  We also examine the
empirical distributions in Figure \ref{f:per1d_hist}, compared to the
density of the QSD (which can be precisely computed by solving an
eigenvalue problem using the formula~\eqref{e:qsd_eig} in this simple
1D situation).  {By $t=1.0$, the qualitative features of
  the distribution are good, and the Gelman-Rubin statistics are less
  than 1.1 for the three observables.  We will see that that $\tol =
  0.1$ gives reasonable results in many cases.}

From this simple experiment, we first observe that the Gelman-Rubin
statistics seems to yield sensible results for assessing the
convergence of the Fleming-Viot particle process. Further inspection
of Figures \ref{f:per1d_grobs} and \ref{f:per1d_hist} suggests
\eqref{e:termination} may be conservative: if the tolerance is too
small, the convergence time may be overestimated compared to what can
be observed on the empirical distribution.

\subsection{Double Well Potential in 2D}
\label{sec:DP2D}

\begin{figure}
  \includegraphics[width=6cm]{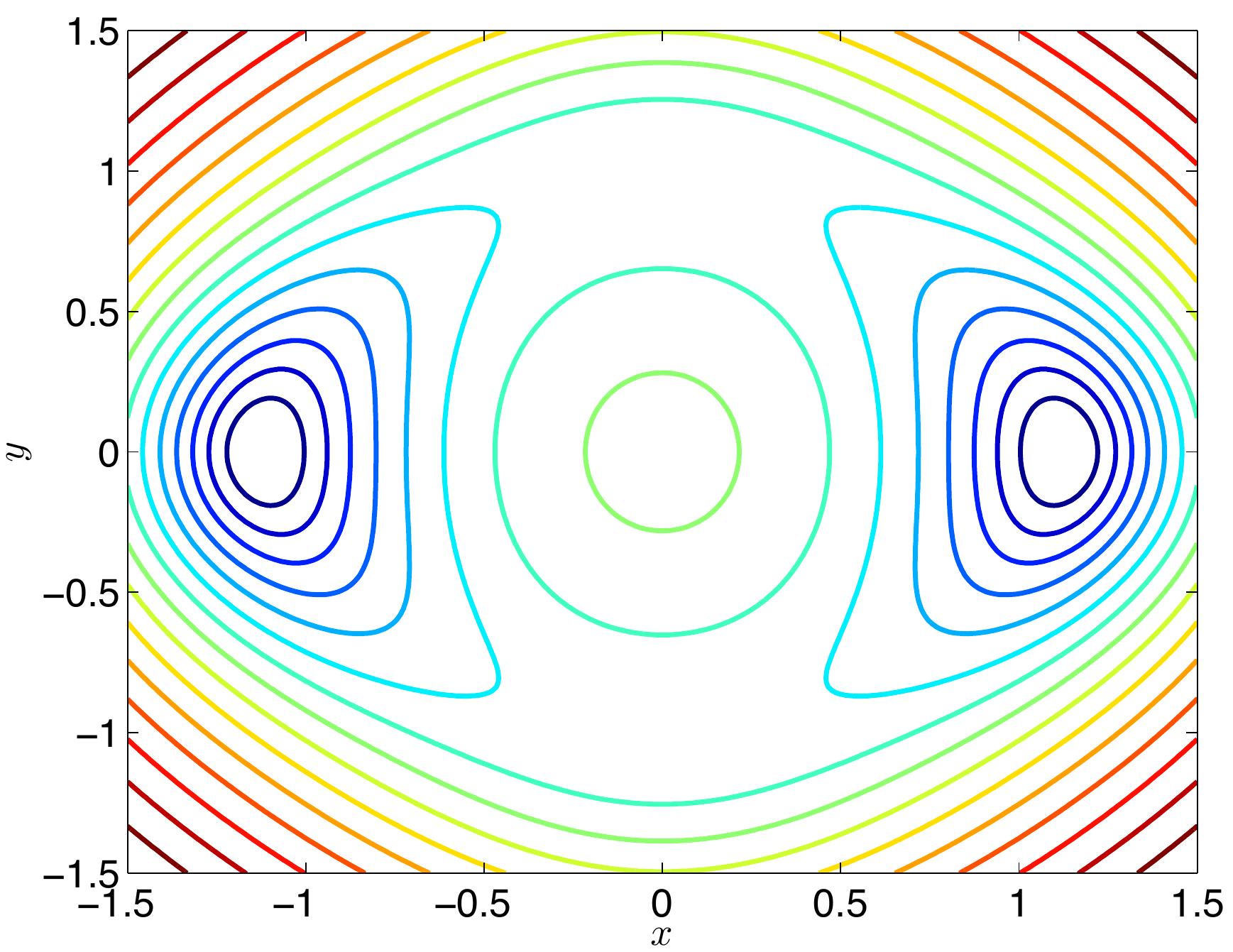}
  \caption{Contours of potential \eqref{e:abf2d_potential}.  Note the
    two minima, near $(\pm 1,0)$.}
  \label{f:abf2d}
\end{figure}

As a second example, we consider the potential
\begin{equation}
  \label{e:abf2d_potential}
  V(x,y) = \tfrac{1}{6}\bracket{4 (1-x^2-y^2)^2 + 2(x^2 -2)^2 + ((x+y)^2
    -1)^2 + ((x-y)^2-1)^2},
\end{equation}
plotted in Figure \ref{f:abf2d}.  Notice that there are two minima,
near $(\pm 1,0)$, along with an internal barrier, centered at the
origin, separating them.  Thus, there are two channels joining the two
minima, with a saddle point in each of these channels. We study this
problem at inverse temperature $\beta = 5$. The aim of this example is
to illustrate the possibility of pseudo-convergence when using
convergence diagnostics, even when it is applied to independent
replicas. We thus concentrate on the sampling of the canonical
distribution with density $Z^{-1}e^{-\beta V(x,y)}$, using independent
realizations, instead of the sampling of the QSD using the
Fleming-Viot particle process.

The observables used in this problem are:
\begin{equation}
  \label{e:abf2d_obs}
  x, \, y, \, V(x,y), \, \norm{ {\bf x}}_{\ell^2}.
\end{equation}

\begin{figure}
  \includegraphics[width=8cm]{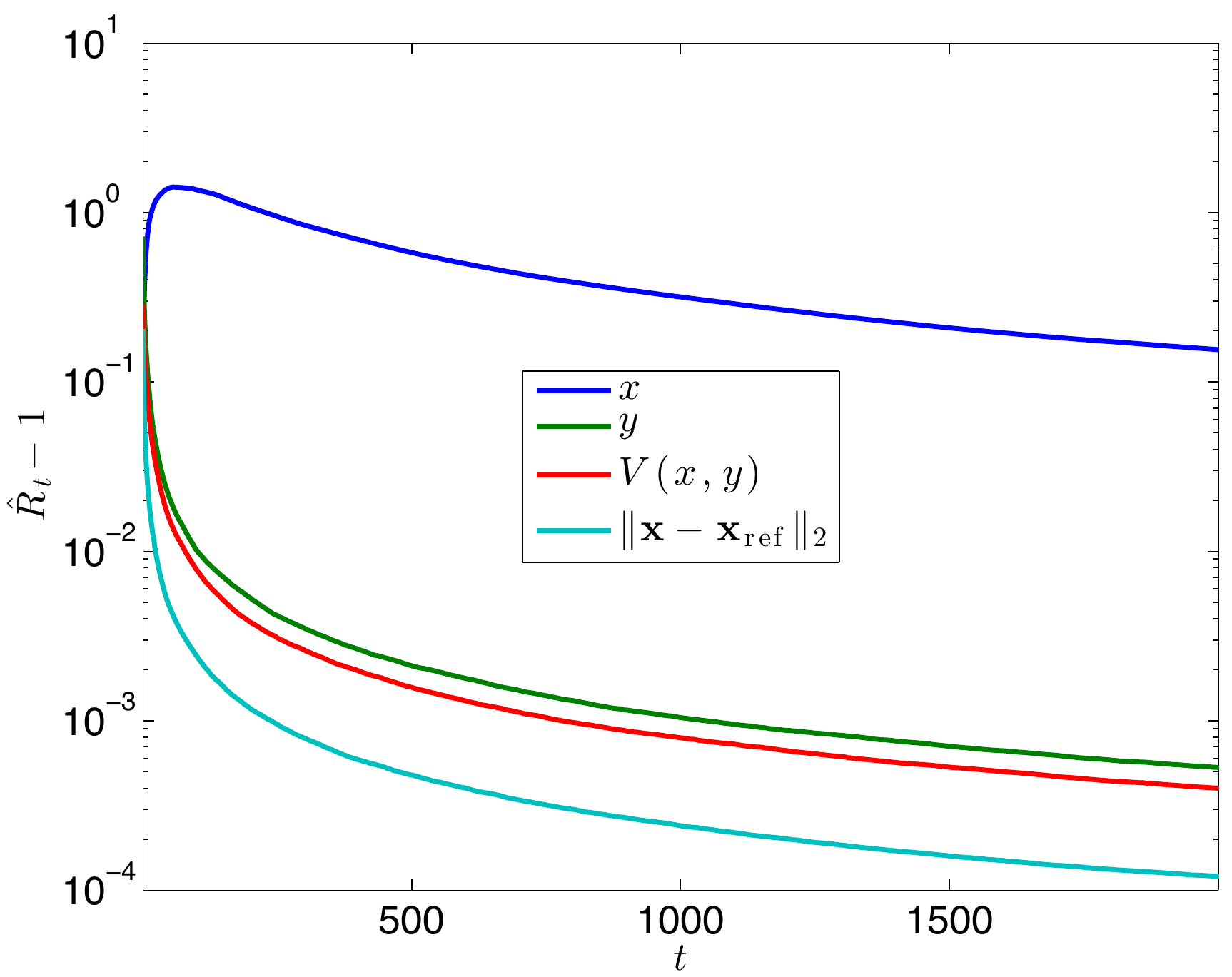}
  \caption{The Gelman-Rubin statistics as a function of time for the
    observables~\eqref{e:abf2d_obs} for $N=10^4$ independent replicas
    following the dynamics~\eqref{e:odlang}. The potential is given
    by~\eqref{e:abf2d_potential}.  Since the replicas begin in the
    left mode, the internal barrier makes it challenging to adequately
    sample the $x$ observable.}
  \label{f:abf2d_grobs}
\end{figure}

\begin{figure}
  \subfigure{\includegraphics[width=6cm]{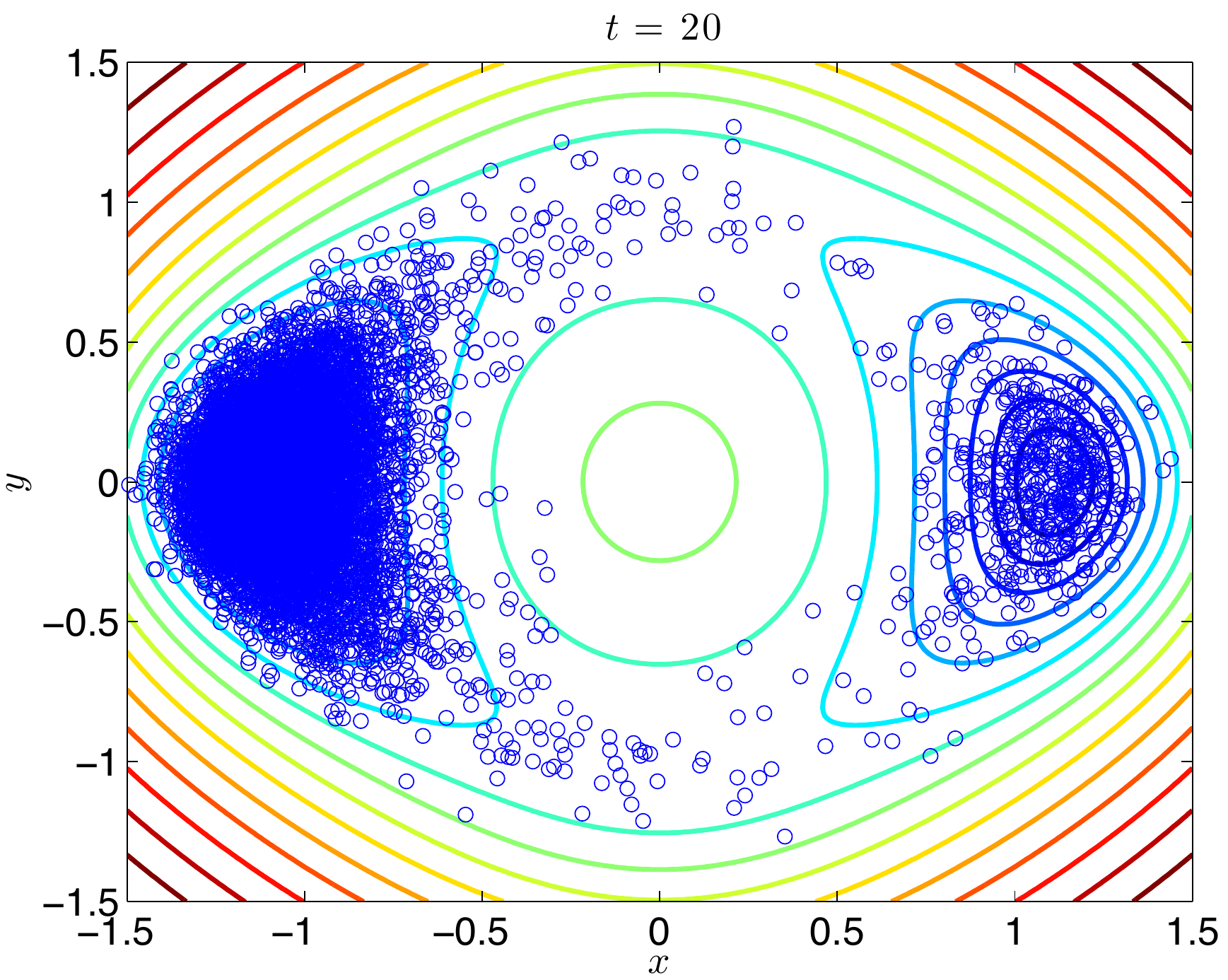}}
  \subfigure{\includegraphics[width=6cm]{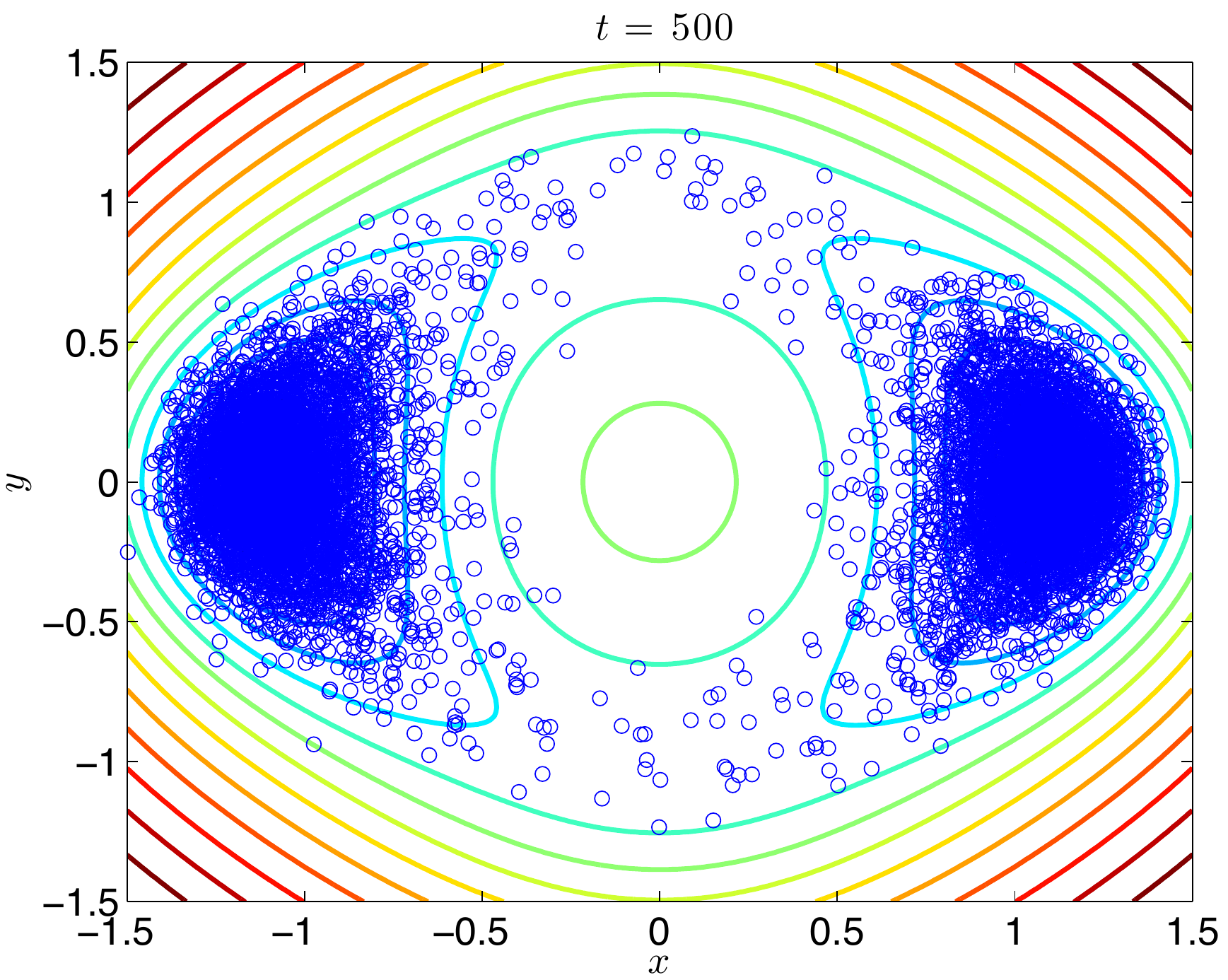}}
  \caption{Empirical distributions from $N=10^4$ independent replicas
    following the dynamics~\eqref{e:odlang} for the
    potential~\eqref{e:abf2d_potential}, see also
    Figure~\ref{f:abf2d_grobs}.}
  \label{f:abf2d_dist}
\end{figure}

Starting our trajectories at $X_0 = (-1,0)$, the Gelman-Rubin
statistics as a function of time appear in Figure~\ref{f:abf2d_grobs},
and the empirical distributions at two specific times are shown in
Figure \ref{f:abf2d_dist}.  We make the following remarks.

First, were we to have neglected the $x$ observable, the Gelman-Rubin
statistics of the other observables would have fallen below $1.1$ by
$t=20$.  But as we can see in Figure \ref{f:abf2d_dist}, this is
completely inadequate for sampling $Z^{-1}e^{-\beta V(x,y)}$. Thus, if
the tolerance is set to $0.1$, the convergence criterion may be
fulfilled before actual convergence if the observables are poorly
chosen. This is a characteristic problem of convergence diagnostics;
they are necessary, but not sufficient to assess convergence.

Let us now consider the $x$ observable, which is sensitive to the
internal barrier. From Figure~\ref{f:abf2d_grobs}, we see that it is
not monotonic and that, even after running till $t=2000$, the
associated statistic still exceeds $1.1$. On the other hand, if we
consider the ensemble at $t=500$, the empirical distribution appears
to equally sample both modes.  Indeed 48\% of the replicas are in the
right basin.  Despite this, the Gelman-Rubin statistic for $x$ is
still relatively large.  This is due to the conservative nature of
\eqref{e:R2}, already mentioned in Section \ref{sec:per1d}.  Once the
ensemble has an appreciable number of samples in each mode, it will
only reach one after all the trajectories have adequately sampled both
modes.

{\subsection{Distribution of $\tcorr$ and $\tphase$}
  \label{s:tdist}

  In the modified ParRep algorithm, $\tcorr=\tphase$, determined by
  \eqref{e:termination}, is a random variable due to $N$ being finite.
  This introduces some uncertainty in the decorrelation time; it may
  be that, for a particular realization, \eqref{e:termination} is
  satisfied at a time which is too small for the replicas, or the
  reference process, to have sampled the QSD.  This can be mitigated
  by using more replicas.  As $N$ tends to infinity, the empirical
  distribution of the Fleming-Viot particle process will converge to
  the law of the process, conditioned on non-extinction, $\mu_t$, as
  in \eqref{e:fvconv1}.  Hence, $\hat{R}_t(\mathcal{O}_j)$ tends to a
  deterministic evolution, and $\tphase$ will approach a deterministic
  value, for the given observables and tolerance.

  To illustrate this, we explored the 1D periodic problem of Section
  \ref{sec:per1d} with $N=10^1, 10^2, 10^3,$ and $10^4$ replicas at
  $\tol = 0.1$.  We ran the process until it satisfied
  \eqref{e:termination}, recording the value of $\tphase$, in each of
  $10^4$ independent realizations, for each of the values of $N$.  The
  results, shown in Figure \ref{f:tqsd_dist}, show the reduction of
  variability as $N$ becomes large.

  In the rest of our examples, appearing in Section
  \ref{s:parrep_examples}, we use $N=100$ replicas.  The typical
  variations of $\tcorr=\tphase$, measured as the ratio of the
  standard deviation to the mean of the samples, are 10\%--20\% of the
  mean.  Both the mean and variance are reported in our data tables,
  below.

  \begin{figure}
    \includegraphics[width=8cm]{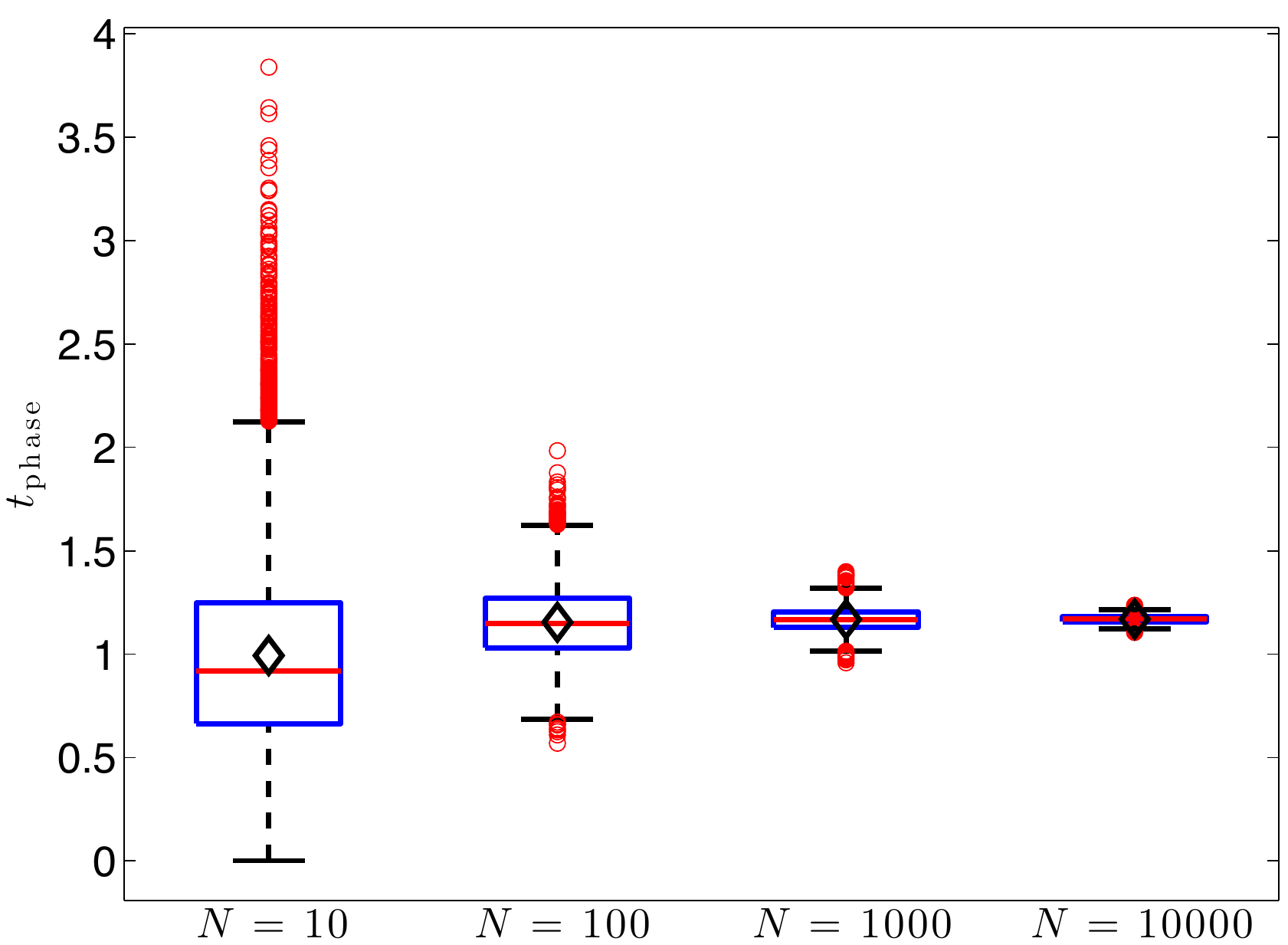}
    \caption{Box plots of $\tphase$ for the 1D periodic problem at
      different values of $N$ with $\tol = 0.1$.  Each distribution
      includes $10^4$ samples.  The box encompasses the middle 50\% of
      the data, the median is the line across the box, and the diamond
      is the mean.  The stems extend out 1.5 times the Interquartile
      Range (the length of the box).  Outliers, outside the stems, are
      shown with circles.}
    \label{f:tqsd_dist}
  \end{figure}}

\subsection{Remarks on Numerical Experiments}

The problems presented above, though simple, demonstrate both the
effectiveness and the caveats of the convergence diagnostics.

First, the convergence diagnostics can be conservative, which is
computationally wasteful. Second, there is the possibility of
pseudo-convergence: as is the case with all convergence diagnostics,
they cannot guarantee stationarity.  Thus, some amount of heuristic
familiarity with the underlying problem is essential to obtain
reasonable results.  One should be careful when choosing the
observables for~\eqref{e:R2}, selecting degrees of freedom which are
associated with the metastable features of the dynamics, as revealed
by the $x$ observable in the preceding example.  Again, the burden is
on the practitioner to be familiar with the system and have some sense
of the relevant observables to test for stationarity.  However, to
set, {\it a priori}, suitable values of $\tcorr$ and $\tphase$
requires more precise knowledge of the system than is needed to apply
convergence diagnostics.  In some sense, this is comparable to the
relationship between {\it a priori} and {\it a posteriori} estimates
in other fields of numerical analysis. {\it a priori} estimates are
generally insufficient for assessing convergence, since they require
some knowledge of the (unknown) solution, such as its regularity.
This is why {\it a posteriori} estimates have been sought out, to
provide error bounds based on computable quantities.

{Hence, the reader may wonder why we did not make use of
  the {\em a priori} convergence estimate \eqref{e:TVconvergence} to
  determine the convergence time to the QSD. The reasons are that we
  lack estimates of both the exponential rate of convergence,
  $(\lambda_2 - \lambda_1)$, and the prefactor, $C(\mu_0)$. This is a
  general challenge to the application of many Markov Chain Monte
  Carlo algorithms: under very weak assumptions, it can be shown that
  the law of a Markov Chain converges to equilibrium exponentially
  fast in time (see, for example, \cite{Meyn:2009aa}).  However, such
  results cannot be used to estimate the convergence, or {\it
    burn-in}, time since the error bound involves difficult to
  estimate quantities.  The need to identify such a time has motivated
  the broad investigation of convergence diagnostics for stochastic
  algorithms.}

\section{Modified ParRep Examples}
\label{s:parrep_examples}

In this section we present a number of numerical results obtained with
the modified ParRep algorithm.  Before presenting the examples, we
review the numerical methods and the parameters used in the
experiments.

\subsection{Common Parameters and Methods}

For each problem we compare the direct, serial simulation to our
proposed algorithm combining ParRep with the Fleming-Viot particle
process and convergence diagnostics, with several different values of
$\tol$ for the stationarity criterion~\eqref{e:termination}.  For each
case of each problem, we perform $10^5$ independent realizations of
the experiment to ensure we have adequate data from which to make
statistical comparisons.  For these ensembles of experiments, which
are performed in parallel, SPRNG~2.0, with the linear congruential
generator, is used as the pseudo-random number generator,
\cite{Mascagni:2000ba}.  In all examples, $N=100$ replicas are used.

The stochastic differential equation \eqref{e:odlang} is discretized
using Euler-Maruyama with a time step of $\Delta t = 10^{-4}$ in all
cases except for the experiments in Section \ref{s:lj7}, where, due to
the computational complexity, $\Delta t = 10^{-3}$.  To minimize
numerical distortion of the exit times due to the time step
discretization, we make use of a correction to time discretized ParRep
presented in~\cite{Aristoff:2014tm}.

Exit distributions obtained from our modified ParRep algorithm and
from the direct serial runs are compared using the two sample
nonparametric {\it Kolmogorov-Smirnov test} (K.-S. Test),
\cite{Gibbons:1997aa}.  Working at the $\alpha =0.05$ level of
significance, this test allows us to determine whether or not we can
accept the null hypothesis, which, in this case, is that the samples
of exit time distributions arose from the same underlying
distribution.  Though the test was formulated for continuous
distributions, it can also be applied to discrete ones, where it will
be conservative with respect to Type 1 statistical errors,
\cite{Crutcher:1975wz,Horn:1977ua}.  In our tables, we report PASS for
instances where the null hypothesis is not rejected (the test cannot
distinguish the ParRep results from the direct serial results), and
FAIL for cases where the null hypothesis is rejected. We also report
the $p$-values, measuring the probability of observing a test
statistic at least as extreme as that which is computed.

In the tables below, we  also report the mean dephasing time
($\mean{t_{\phase}}$), its variance ($\Var(t_{\phase})$), the average
exit time ($\mean{T}$), and the percentage of realizations for which
the decorrelation/dephasing step were successful ($\%$ Dephased). The
average and variance of $t_\phase$ are computed only over
realizations for which the reference process decorrelates.

{Finally, we report the speedup of the algorithm, which we
  define as
  \begin{equation}
    \label{e:speedup}
    \text{Speedup} \equiv\frac{\text{Physical Exit
        Time}}{\text{Computational Time}}.
  \end{equation}
  The computational time is either the exit time of the reference
  process (if the exit occurs during the decorrelation step), or the
  sum of the time spent in the decorrelation/dephasing step and in the
  parallel step (if the exit occurs during the parallel
  step). Formula~\eqref{e:speedup} measures the speedup of the
  algorithm over direct numerical simulation of the exit event.  Let
  us heuristically consider the mean behavior of \eqref{e:speedup}.
  Let $p$ denote the probability that the system decorrelates, so that
  $1-p$ corresponds to the probability that the reference process
  exits before \eqref{e:termination} is satisfied.  For realizations
  where the reference process exited before decorrelating,  the speedup is unity.  For
  realizations which exit in the parallel step, let $\mean{\tphase}$
  and $\mean{\tpar}$ denote the average dephasing and parallel exit
  times.  Roughly, these values relate to the eigenvalues, introduced
  in \eqref{e:spectrum}, and the choice of $N$ as
  \begin{equation}
    \mean{\tphase}\approx \frac{C_\phase}{\lambda_2 -\lambda_1},
    \quad \mean{\tpar}\approx \frac{1}{N\lambda_1}
  \end{equation}
  To obtain $\mean{\tpar}$, we have assumed that the dephased
  distribution is a good approximation of the QSD.  The constant
  $C_\phase$ is determined by the initial distribution, the tolerance,
  the observables, and the temperature of the system.

  Thus, assuming these random variables are not too strongly
  correlated, we have that, roughly, 
  \begin{equation}
    \mean{\text{Speedup}} \approx (1-p) + p \frac{\mean{\tphase} +N
      \mean{\tpar} }{\mean{\tphase} + \mean{\tpar} }\approx (1-p) + p \frac{\frac{C_\phase}{\lambda_2 -\lambda_1} +
      \frac{1}{\lambda_1} }{\frac{C_\phase}{\lambda_2 -\lambda_1} +\frac{1}{N\lambda_1} }.
  \end{equation}
  As $N$ tends to infinity, the maximum speedup monotonically
  saturates:
  \begin{equation}
    \label{e:speedup_limit}
    \lim_{N\to \infty}\mean{\text{Speedup}} \approx 1 + p \frac{1}{C_{\phase}}\frac{\lambda_2-\lambda_1}{\lambda_1}.
  \end{equation}
  Alternatively, if
  \begin{equation}
    \label{e:Nlimit}
    N \ll \frac{1}{C_\phase}\frac{\lambda_2 - \lambda_1}{\lambda_1},
  \end{equation}
  we can expand the expression to recover the linear speedup regime,
  \begin{equation}
    \label{e:linspeedu}
    \mean{\text{Speedup}} \approx (1-p) + p N \paren{1 +C_\phase \frac{
        \lambda_1}{\lambda_2 - \lambda_1}}.
  \end{equation} 
  For strongly metastable systems, we expect that $p$ will be close to
  one (the reference process almost always decorrelates) and that the relative
  spectral gap will be large, $\lambda_2 - \lambda_1 \gg \lambda_1$
  (the typical exit time is much larger than the time to converge to
  the QSD). The constant $C_\phase$ will also depend on the degree of
  metastability, though the dependence is less clear. Thus, we expect there will be a large range of $N$ satisfying
  \eqref{e:Nlimit}, and a linear speedup will be observed.  However,
  there is also an intrinsic limitation on the performance of ParRep
  that depends on the strength of metastability of each state.

  The ratio in \eqref{e:speedup_limit} has a natural interpretation.
  As $N$ becomes large in ParRep, the cost of finding the first exit
  during the parallel step is driven to zero, and what remains is the
  cost of decorrelation and depahsing.  These pieces of the algorithm
  have a time scale of $(\lambda_2 -\lambda_1)^{-1}$.  For a direct
  serial simulation, in a sufficiently metastable system, the cost is
  to find the first exit, which has a time scale of $\lambda_1^{-1}$.
}

Since dephasing and decorrelation terminate at a finite time, there
will always be some bias in the exit distributions, which is in
general difficult to distinguish from the statistical noise.  When the
number of independent realizations of the experiment becomes
sufficiently large, the statistical tests may detect this bias and
identify two distinct distributions, even if the difference between
the two distributions is tiny.  Two things are thus to be expected
when comparing the exit time distributions obtained by ParRep and
direct numerical simulations:
\begin{itemize}

\item {For a fixed number of independent realizations, as
    the tolerance is reduced to zero, the $p$-value of the K.-S.  test
    goes to $1$, and the null hypothesis is accepted; the
    distributions will be the same.  This is because, as the tolerance
    becomes more stringent, it takes longer to reach stationarity.  On
    the one hand, there is then a higher chance that the reference
    process exits during the joint decorrelation/dephasing step.  On
    the other, for the realizations which manage to satisfy
    convergence criterion~\eqref{e:termination}, the distribution will
    be closer to that of the QSD.  Both these effects induce ParRep to
    better replicate the true, unaccelerated distributions.}

\item For fixed tolerance, as the number of realizations tends to
  $+\infty$, the $p$-value of the K.-S.  tends to zero, and the null
  hypothesis is rejected. This is because as the statistical noise is
  reduced, the bias due to a finite stationarity time for a given
  tolerance becomes apparent.  This is born out in numerical
  experiments.

\end{itemize}
In addition, due to sampling variability, the $p$-value may not vary
monotonically as a function of $\tol$ or of the number of realizations

\subsection{Periodic Potential in 2D}
\label{s:per2d}

\begin{figure}
  \includegraphics[width=8cm]{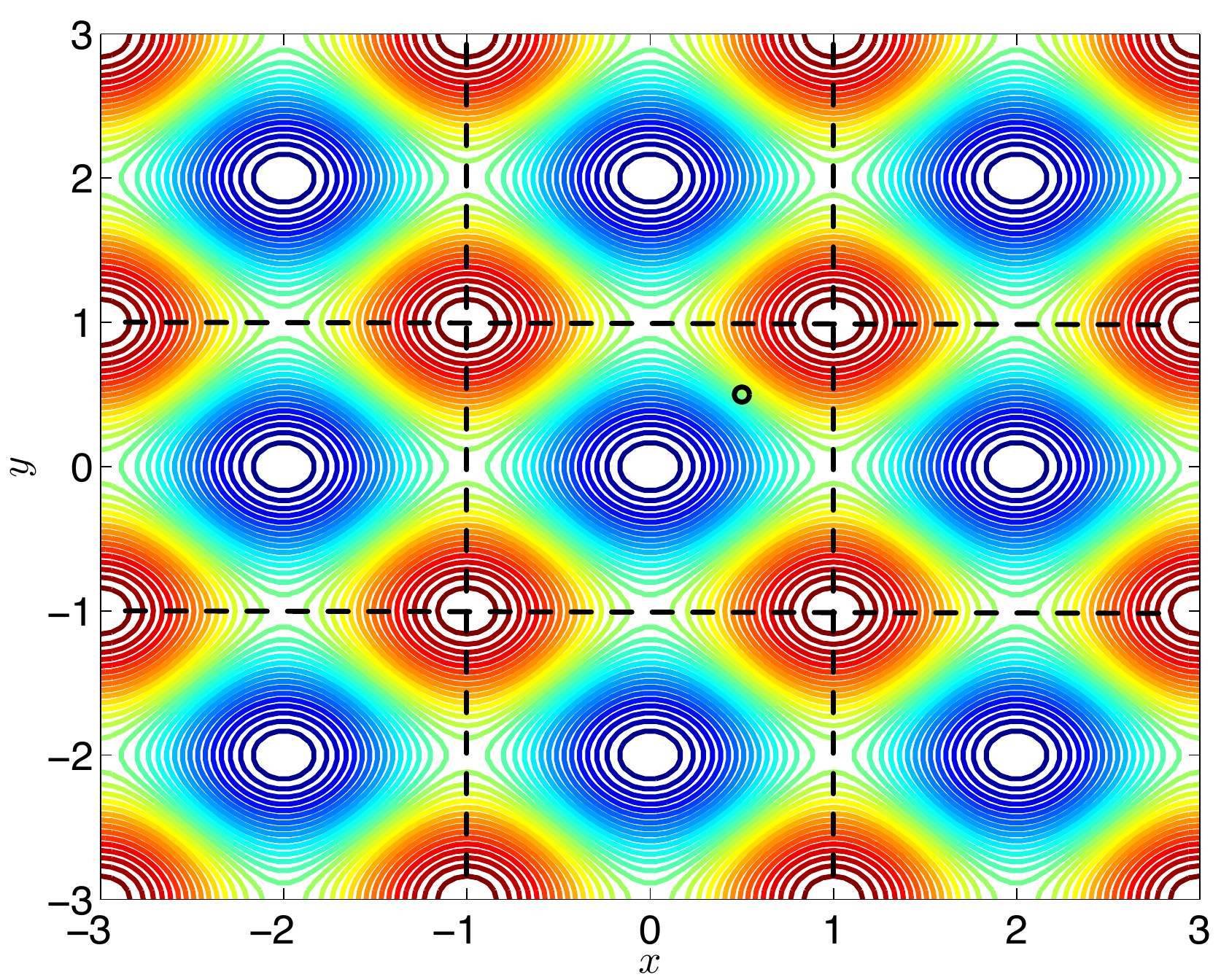}
  \caption{The partition of the configuration space into states for
    the potential~\eqref{e:per2d_potential}. The circle $\circ$
    indicates the initial condition, $(.5,.5)$.}
  \label{f:per2d_multi}
\end{figure}

As a first example, we consider the two dimensional periodic potential
\begin{equation}
  \label{e:per2d_potential}
  V(x,y) = - \cos(\pi x) - \cos(\pi y),
\end{equation}
simulated at $\beta = 3$.  The states are defined as the translates of
$(-1,1)^2$, seen in Figure \ref{f:per2d_multi}.  The observables used
with the Gelman-Rubin statistics in~\eqref{e:R2} are
\begin{equation}
  \label{e:per2d_obs}
  x, \quad y, \quad V(x,y), \quad\norm{ {\bf x}- {\bf x}_{\refe} }_{\ell^2},
\end{equation}
where ${\bf x}_{\refe}$ will be the local minimum of the state under
consideration.

\subsubsection{Escape from a Single State}

\begin{figure}
  \subfigure{\includegraphics[width=7cm]{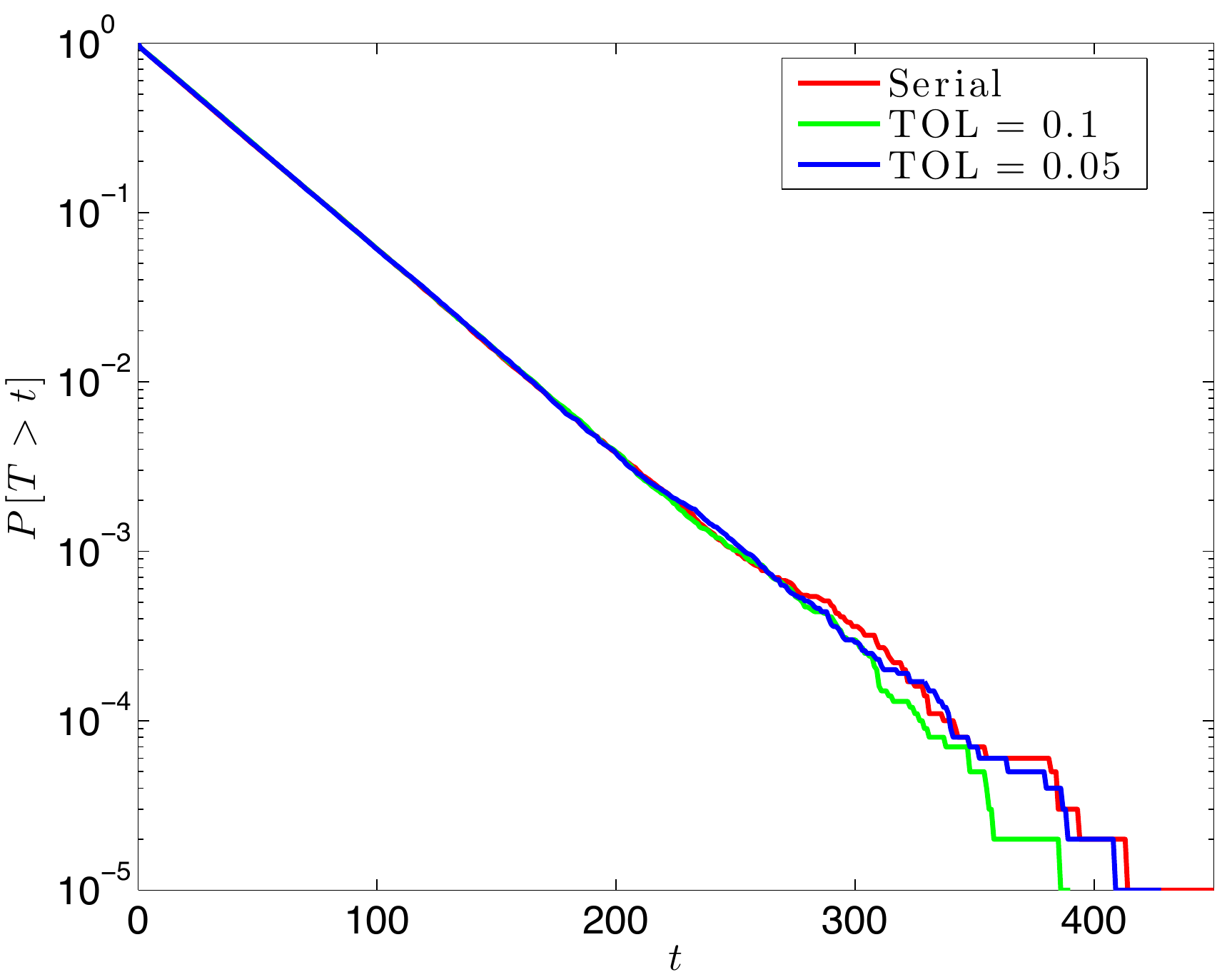}}
  \subfigure{\includegraphics[width=7cm]{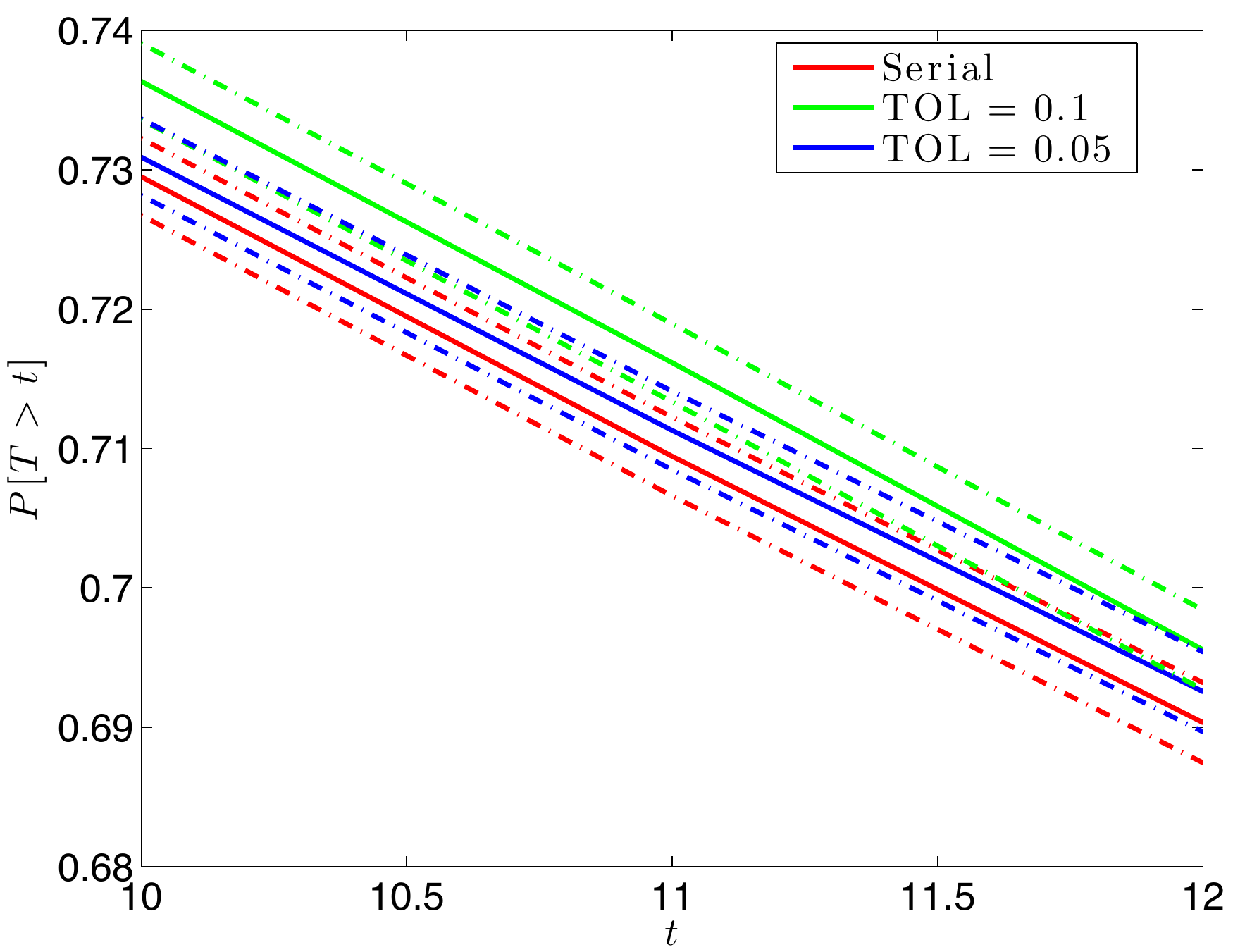}}
  \subfigure{\includegraphics[width=7cm]{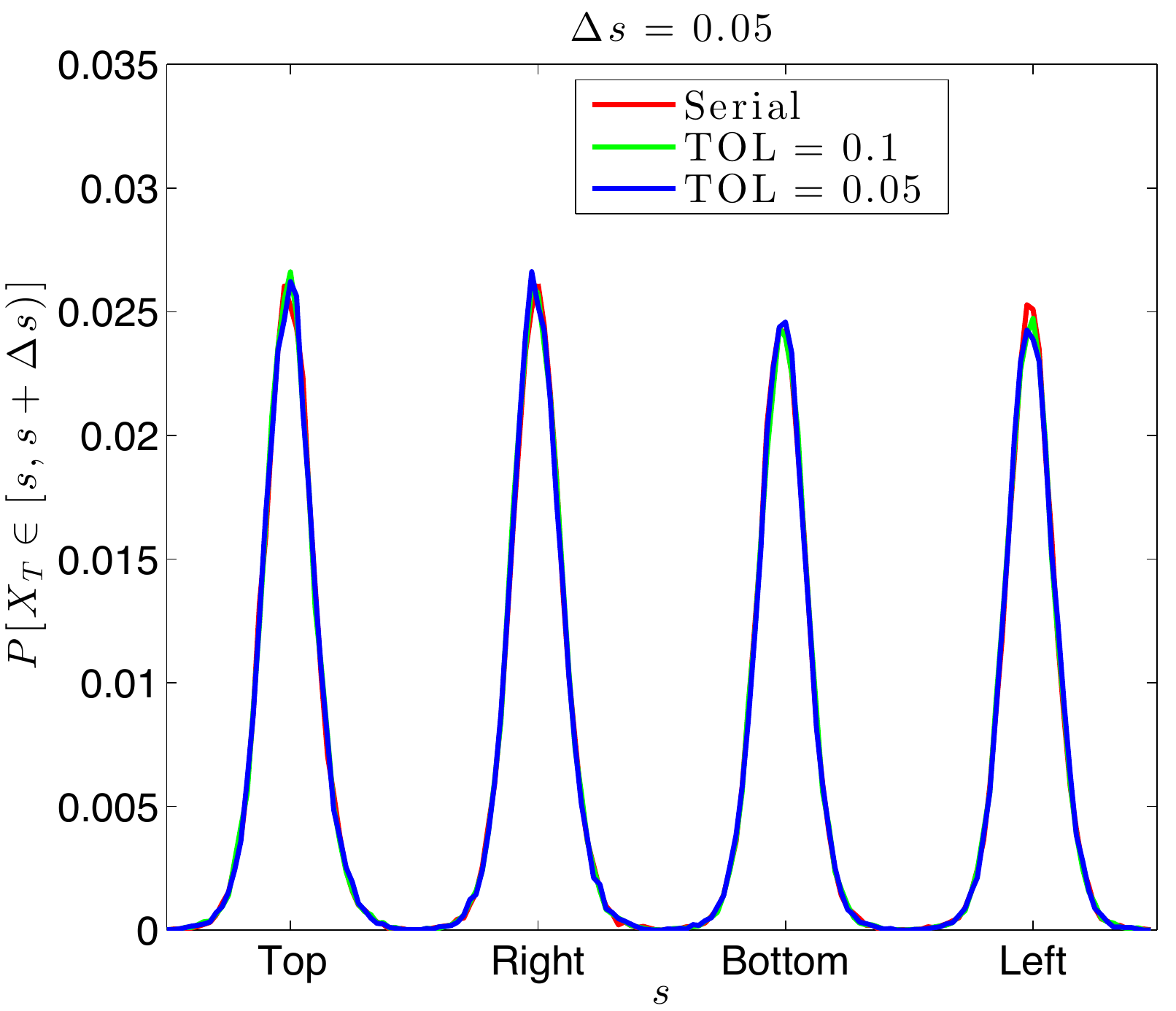}}
  \caption{{\bf Periodic Potential in 2D--Single Escape:} Exit
    distributions from $\Omega = (-1,1)^2$ for the 2D Periodic
    potential \eqref{e:per2d_potential}.  On the right is a
    magnification of the figure from the left; notice that the
    confidence intervals do not overlap at $\tol = .1$.  For the
    hitting point distribution, the boundary of $\Omega$ has been
    treated as a 1D manifold parametrized clockwise, with the labels
    ``Top'', ``Right'', ``Bottom'', and ``Left'' corresponding
    respectively to the edges $ (-1,1)\times \{1\}$, $\{1\}\times
    (-1,1)$, {\it etc.}}
  \label{f:per2d_exit}
\end{figure}

We first consider the problem of escaping from the state $\Omega =
(-1,1)^2$ with initial condition $X_0 = (0.5,0.5)$. We run $10^5$
independent realizations of Serial and ParRep simulation, with
quantitative data appearing in Table~\ref{t:per2d}, and the
distributions in Figures~\ref{f:per2d_exit}.

At the tolerance value $\tol=0.2$, there is already good agreement in
terms of first hitting point distributions, as seen in Table
\ref{e:per2d_potential}.  To compare the first hitting point
distributions, the boundary of $(-1,1)^2$ is treated as a one
dimensional manifold composed of the four edges (top, right, bottom,
and left); this is the third plot of Figure \ref{f:per2d_exit}.  We
indeed have good agreement.

Unlike the exit point distribution, the exit time distribution seems
to require a more stringent tolerance. There are statistical errors at
$\tol =0.2$, which are mitigated as the tolerance is reduced. At $\tol
=0.05$, the K.-S. test shows agreement. However, based on the first
two of Figure \ref{f:per2d_exit}, the exit time distribution is in
good qualitative agreement even for $\tol =0.1$, in the sense that the
$L^1$ distance between the two distributions is small and the decay
rates are in agreement.

Notice that for $\tol =0 .05$, there is good performance in the sense
that 84\% of the realizations exited during the parallel step, and the
statistical agreement is high.  The corresponding speedup is 6.25.

With regard to our earlier comment that, at fixed tolerance, as the
sample size increases, the $p$-values tend to decrease, we consider
the case of $\tol = 0.1$ here.  If we compare the first $10^2$,
$10^3$, $10^4$ and $10^5$ experimental realizations against the
corresponding serial results, we obtain $p$-values of 0.68, 0.18,
0.0099, and 0.0037.

\begin{table}
  \caption{{\bf Periodic Potential in 2D--Single Escape:} Comparison of ParRep and an unaccelerated serial process
    escaping from $(-1,1)^2$ with $X_0=(0.5,0.5)$ for potential
    \eqref{e:per2d_potential}.}
  \begin{tabular}{l l c c c c c l l}
    \hline \hline\\
    Method  & $\tol$ & $\mean{t_{\phase}}$ &$\Var(t_{\phase})$  & $\mean{T}$ &
    $\mean{\text{Speedup}}$ & \% Dephased & $X_T$
    K.-S. Test ($p$) &  $T$ K.-S. Test ($p$)\\
    \hline\\
    Serial   &-- & --& --&34.8& --  & -- & -- & -- \\
    ParRep&0.2 & 1.12& 0.0267&35.3& 20.8  & 93.5\% & PASS (0.59) & FAIL ($2.0 \times 10^{-5}$) \\
    ParRep&0.1 & 2.43 &0.103 & 35.0& 11.6 & 90.2\% & PASS (0.48) & FAIL (0.0037)\\
    ParRep&0.05& 5.10& 0.393&34.8 &6.25 & 83.6\%& PASS (0.52)  & PASS (0.33) \\
    ParRep&0.01 & 26.2& 9.08&34.8 & 1.63 & 46.5\% & PASS (0.27) & PASS (0.42)\\
    \hline
  \end{tabular}
  \label{t:per2d}
\end{table}

\subsubsection{Escape from a Region Containing Multiple States}

\begin{figure}
  \subfigure{\includegraphics[width=7cm]{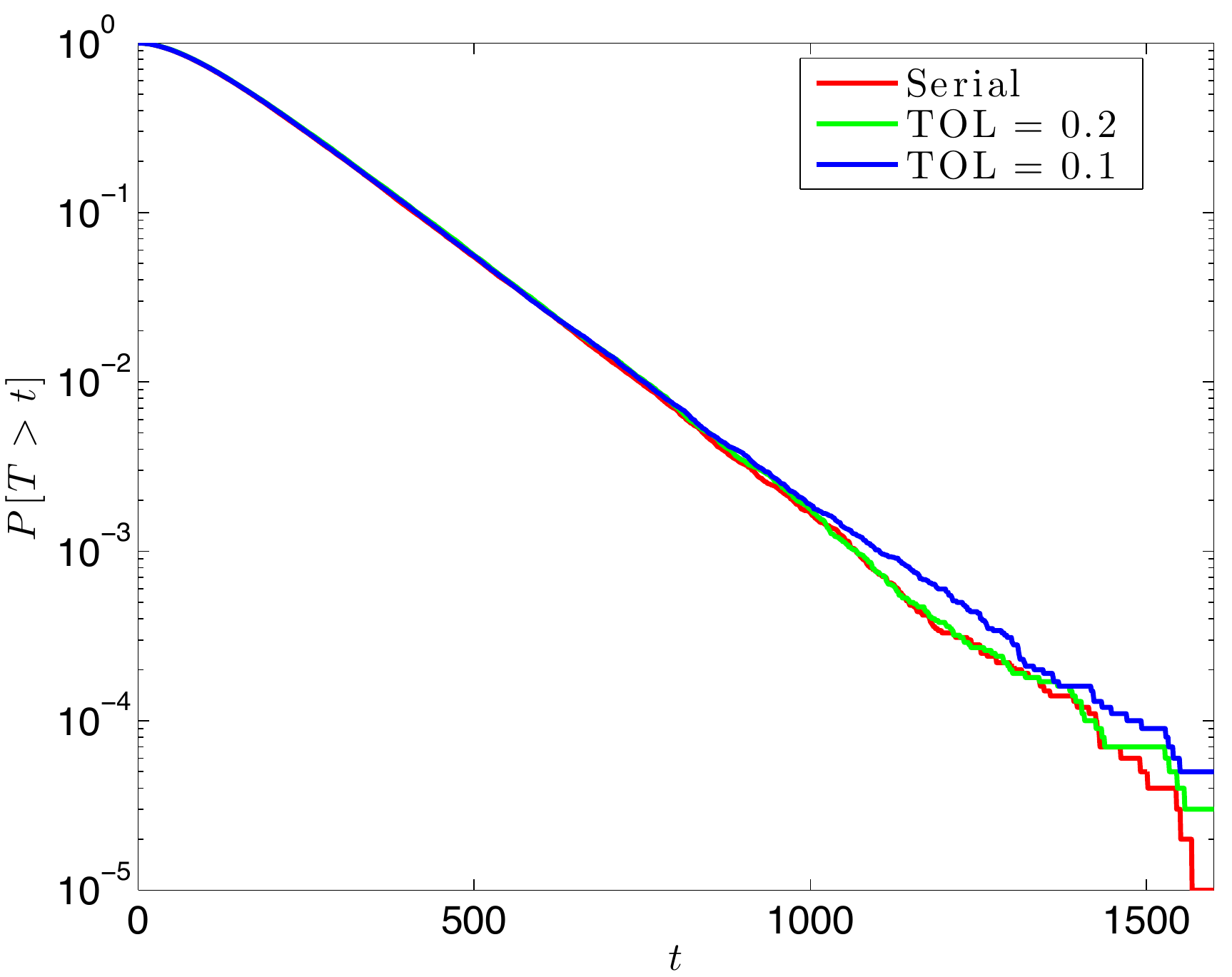}}
  \subfigure{\includegraphics[width=7cm]{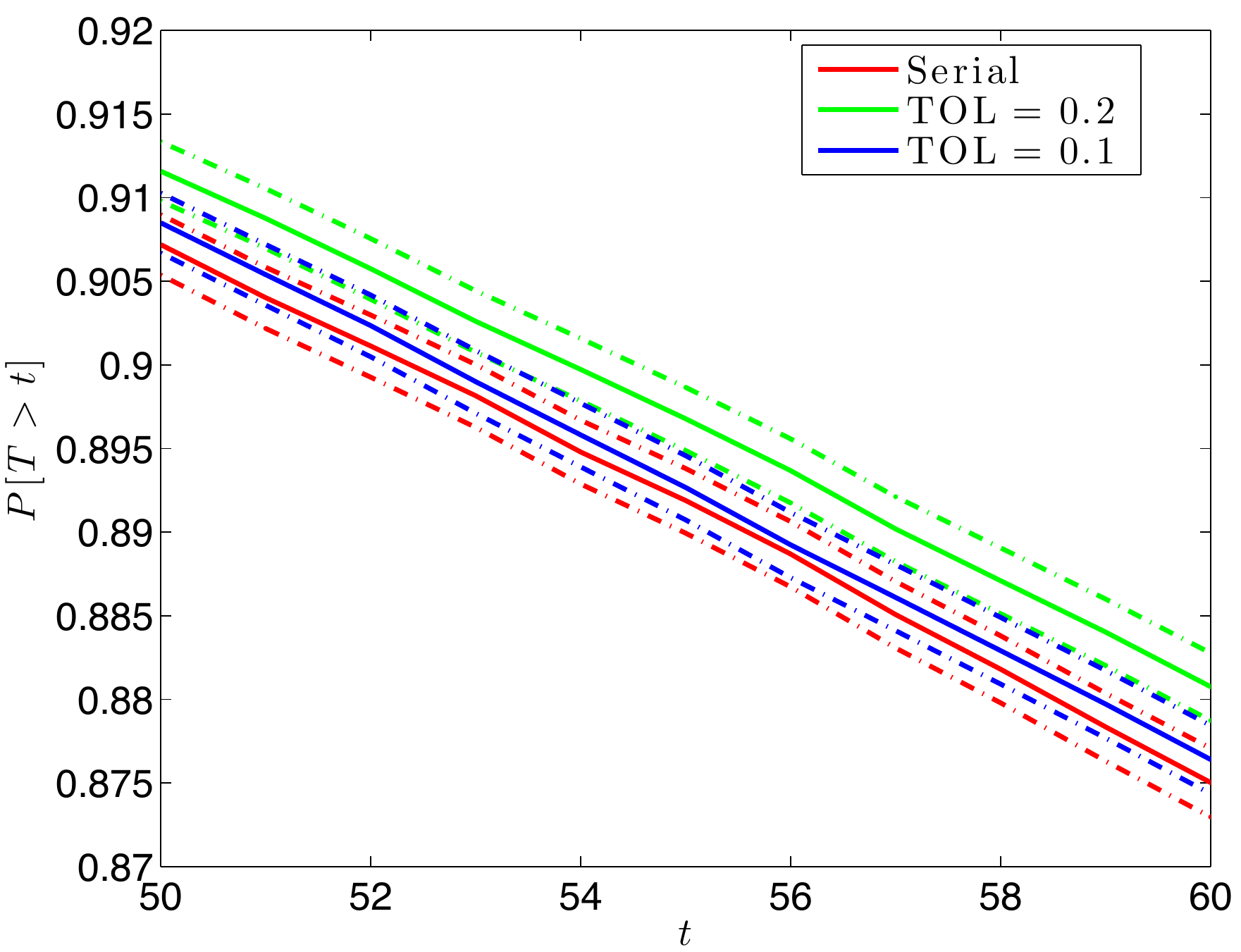}}
  \subfigure{\includegraphics[width=7cm]{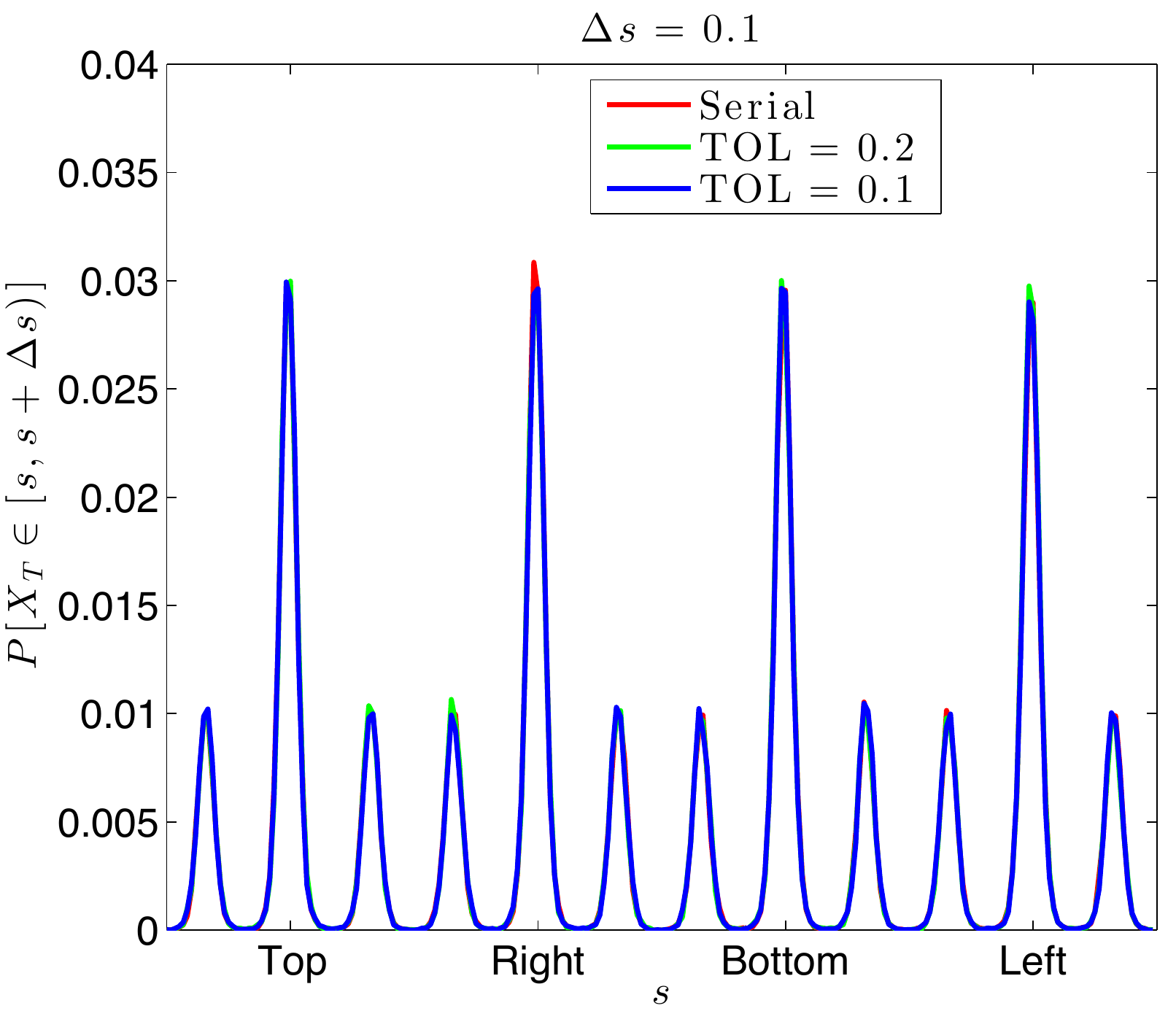}}
  \caption{{\bf Periodic Potential in 2D--Multiple Escapes:} Exit time
    distribution for the 2D periodic potential
    \eqref{e:per2d_potential} when ParRep is used to escape from
    $(-3,3)^2$.  On the right is a magnification of the figure from
    the left; notice that the confidence intervals do not overlap at
    $\tol = .2$. For the hitting point distribution, the boundary of
    $(-3,3)^2$ has been treated as a 1D manifold parametrized
    clockwise, with the labels ``Top'', ``Right'', ``Bottom'', and
    ``Left'' corresponding respectively to the edges $ (-3,3)\times
    \{3\}$, $\{3\}\times (-3,3)$, {\it etc.}}
  \label{f:per2d_multi_exit}
\end{figure}

Next, we consider the problem of escaping from the region $(-3,3)^2$,
running our modified ParRep algorithm over the states indicated in
Figure~\ref{f:per2d_multi}. The results are reported in
Table~\ref{t:per2d_multi} and Figure~\ref{f:per2d_multi_exit}.
Beneath the value $\tol=0.1$, the algorithm is statistically
consistent, both in terms of first hitting point and first exit time
distributions.

\begin{table}
  \caption{{\bf Periodic Potential in 2D--Multiple Escapes:} Comparison of ParRep and an unaccelerated serial process
    escaping from $(-3,3)^2$ with $X_0=(0.5,0.5)$  for potential~\eqref{e:per2d_potential}.}
  \begin{tabular}{ c l l}
    \hline \hline\\
    $\tol$ & $X_T$ K.-S. Test ($p$) &  $T$ K.-S. Test ($p$)\\
    \hline\\
    0.2 & PASS (0.56) & FAIL ($5.2\times 10^{-3}$)\\
    0.1 & PASS (0.51) & PASS (0.24) \\
    0.05 & PASS (0.53) & PASS (0.59) \\
    0.01 & PASS (0.11) &  PASS (0.87) \\
    \hline
  \end{tabular}
  \label{t:per2d_multi}
\end{table}

\subsection{Entropic Barriers in 2D}
\label{s:ent2d}

\begin{figure}
  \includegraphics[width=8cm]{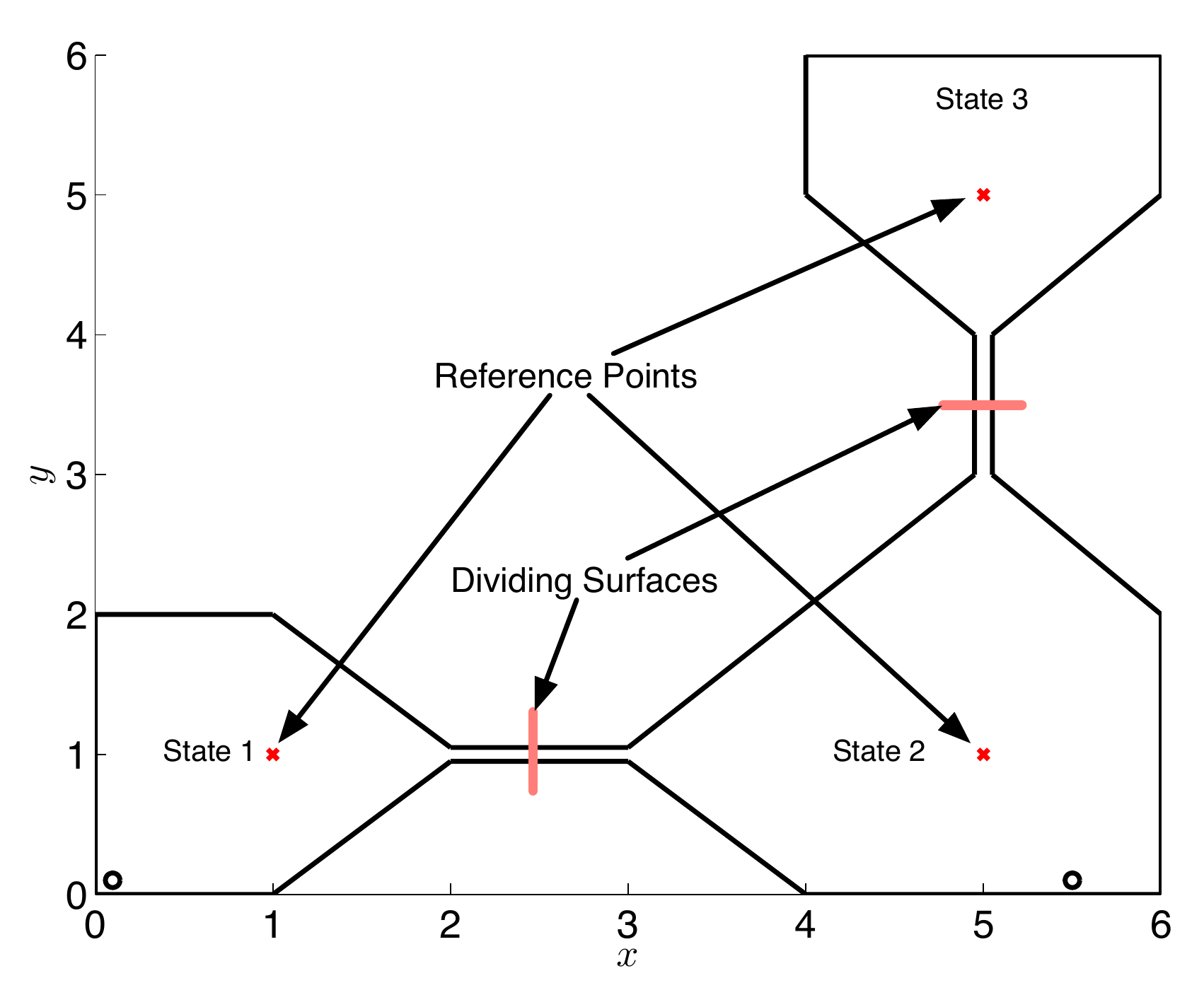}
  \caption{A domain with entropic barriers.  The trajectory begins at
    $(5.5,0.1)$ for the experiment in Section \ref{s:ent_exper1} and
    at $(0.1,0.1)$ for the experiment in Section \ref{s:ent_exper2}
    (see the~$\circ$ symbols).  The width of the necks between states
    is 0.1, and the dividing surfaces between states are indicated by
    the short segments.  A separate reference point is used in each
    state for the observable $\norm{{\bf x} - {\bf x}_{\refe}}$.}
  \label{f:ent2d_domain}
\end{figure}

Let us now consider a test case with entropic barriers.  Consider pure
Brownian motion in the domain represented in Figure
\ref{f:ent2d_domain}, with reflecting boundary conditions.  Here, the
observables are
\begin{equation}
  \label{e:ent2d_obs}
  x, \, y, \, \norm{{\bf x} - {\bf x}_{\refe}}_{\ell^2},
\end{equation}
with the reference points indicated in Figure \ref{f:ent2d_domain}.

\subsubsection{Escape from State 2}
\label{s:ent_exper1}

\begin{figure}
  \subfigure{\includegraphics[width=7cm]{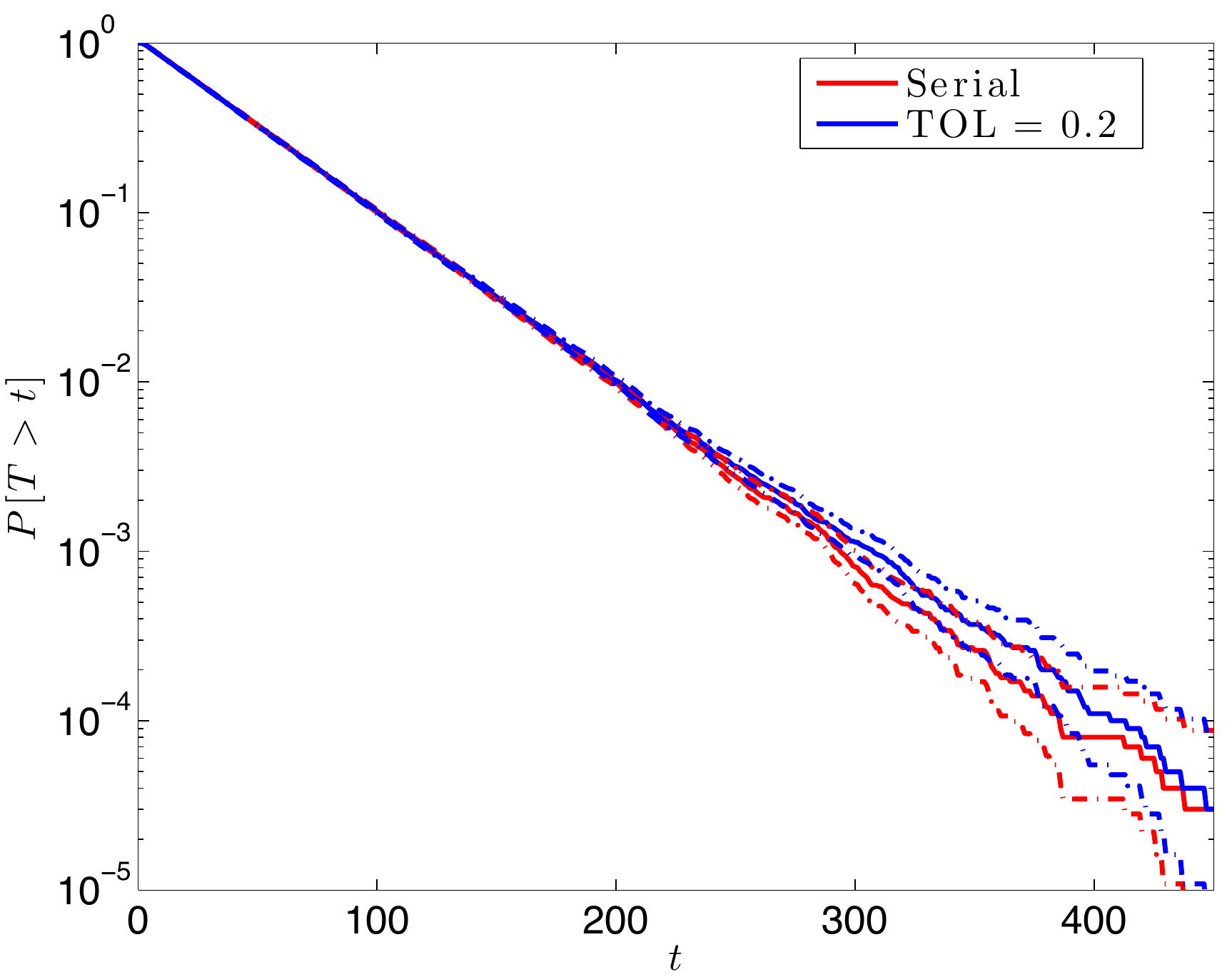}}
  \caption{{\bf Entropic Barrier in 2D--Single Escape:} Exit time
    distribution from state 2 of Figure \ref{f:ent2d_domain} with
    initial conditions $(5.5,.1)$, along with confidence intervals.}
  \label{f:ent2d_single}
\end{figure}

As a first experiment, we look for first escapes from state 2, with
initial condition $(5.5,0.1)$.  Quantitative results appear in Table
\ref{t:ent2d_single} and the exit time distributions are plotted in
Figure \ref{f:ent2d_single}.  For this problem, there is good
statistical agreement even at $\tol = 0.2$.  As the hitting point
distributions across each of the two channels are nearly uniform, we
only report the probability of passing into state 3 versus state 1.
Note that in the $\tol=0.01$ case, the exit is almost always due to
the reference process escaping before stationarity is achieved. This
is an example of setting the parameter so stringently as to render
ParRep inefficient.

\begin{table}
  \caption{{\bf Entropic Barrier in 2D--Single Escape:} Comparison of ParRep and an unaccelerated serial process
    escaping from state 2 of Figure \ref{f:ent2d_domain} with initial
    condition $(5.5,.1)$. }
  \begin{tabular}{l l c c c c c c l}
    \hline \hline\\
    Method  & $\tol$ & $\mean{t_{\phase}}$& $\Var(\tphase)$ & $\mean{T}$&
    $\mean{\text{Speedup}}$ & \% Dephased& 
    $\prob[\text{state 3}]$ &  $T$ K.-S. Test ($p$)\\
    \hline\\
    Serial     & -- &--& --& 45.0 &-- &-- &  $(0.493, 0.500)$ &--\\
    ParRep & 0.2 & 12.7& 2.31& 44.8& 3.46 & 77.9\%& $(0.495, 0.502)$ & PASS (0.90)\\
    ParRep &0.1 & 25.4& 8.75&45.1&  1.95 & 58.0\% & $(0.493, 0.499)$ & PASS (0.51)\\
    ParRep & 0.05& 50.3& 33.3&45.1&  1.27 &32.3\% & $(0.497, 0.503)$ & PASS (0.41)\\
    ParRep & 0.01& 236.0 & 658.0& 45.1& 1.00 & 0.342\%& $(0.496, 0.502)$ & PASS (0.31)\\
    \hline
  \end{tabular}
  \label{t:ent2d_single}
\end{table}

\subsubsection{Getting to State 3 from State 1}
\label{s:ent_exper2}

\begin{figure}
  \subfigure{\includegraphics[width=7cm]{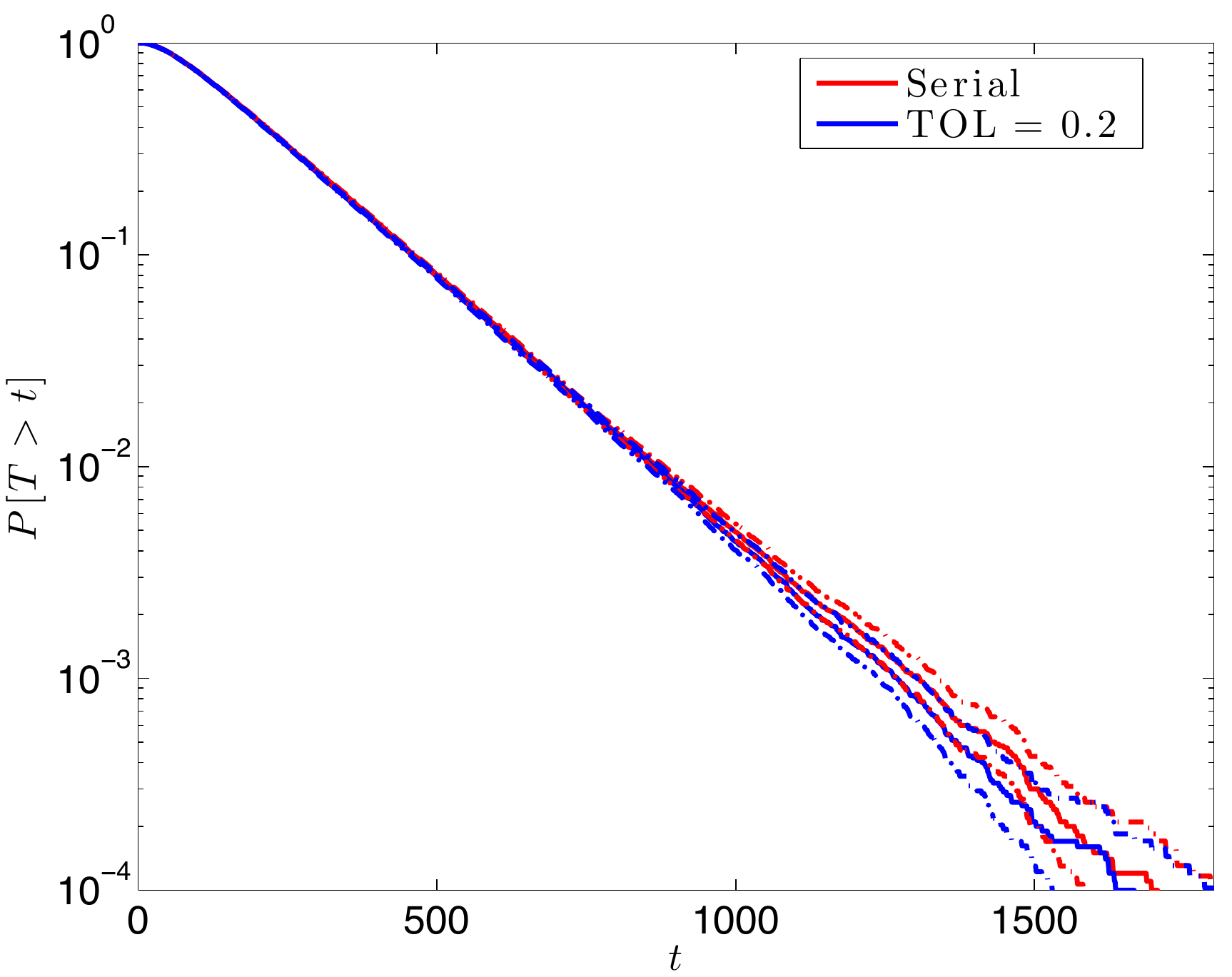}}
  \caption{{\bf Entropic Barrier in 2D--Multiple Escapes:} Exit time
    distributions to get to state 3 from state 1 in the entropic
    barrier problem pictured in Figure \ref{f:ent2d_domain}, along
    with confidence intervals.  The full algorithm is applied in each
    state the trajectory visits.}
  \label{f:ent2d_multi}
\end{figure}

As a second experiment, we start the trajectories in state 1 at
$(0.1,0.1)$, and examine how long it takes to get to state 3, running
the modified ParRep algorithm over the states represented on
Figure~\ref{f:ent2d_domain}.  The results of this experiment are given
in Table \ref{t:ent2d_multi} and Figure \ref{f:ent2d_multi}.  Again,
there is a very good statistical agreement in all cases.

\begin{table}
  \caption{{\bf Entropic Barrier in 2D--Multiple Escapes:} Comparison of ParRep and an unaccelerated serial process
    getting from state 1 to state 3 with initial condition $(.1,.1)$; see
    Figure \ref{f:ent2d_domain}.}
  \label{t:ent2d_multi}
  \begin{tabular}{c l}
    \hline \hline\\
    $\tol$ &$T$ K.-S. Test ($p$)\\
    \hline\\
    0.2 & PASS (0.070)\\
    0.1 & PASS (0.14) \\
    0.05 & PASS (0.090) \\
    0.01 & PASS (0.74) \\
    \hline
  \end{tabular}
\end{table}

\subsection{Lennard-Jones Clusters}
\label{s:lj7}

For a more realistic problem, we consider a Lennard-Jones cluster of
seven atoms in 2D, denoted $\LJ$.  The potential used in this problem
is then:
\begin{equation}
  \label{e:LJ7potential}
  V_7^{\rm 2D}({\bf x}) = \frac{1}{2}\sum_{i\neq j}\phi(\norm{{\bf x}_i
    - {\bf x}_j}), \quad {\bf x}=({\bf x}_1,{\bf x}_2,\ldots, {\bf
    x}_7), \quad \phi(r) = r^{-12} - 2 r^{-6}.
\end{equation}
Due to the computational cost associated with this problem, the time
step size is increased to $\Delta t = 10^{-3}$.  A smaller time step
could have been used, but we would have been precluded from obtaining
large samples of the problem.

\begin{figure}
  \subfigure[$C_0$,
  $V=-12.53$]{\includegraphics[width=3cm]{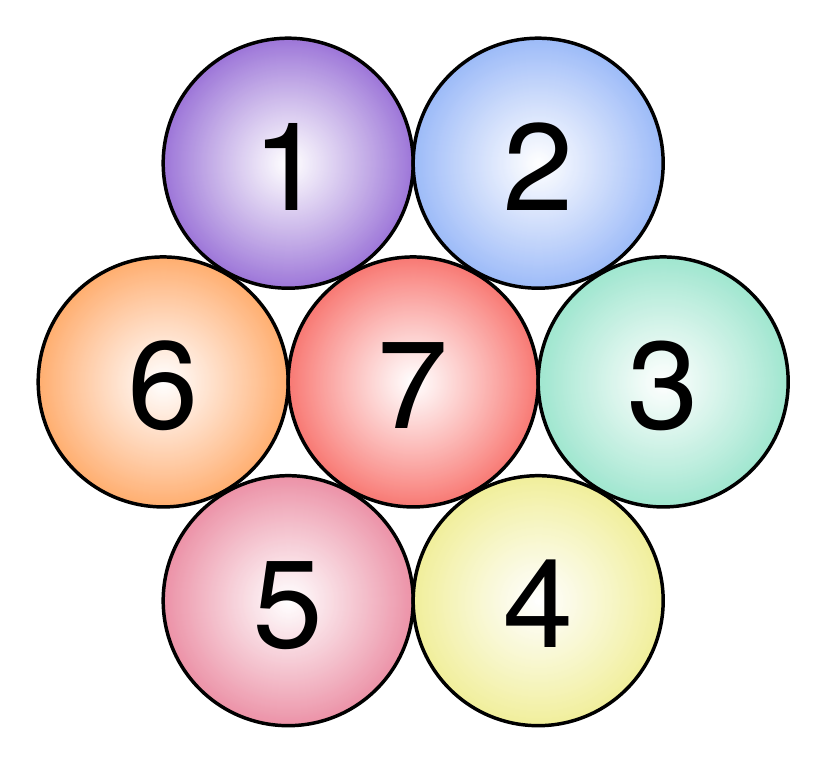}}
  \subfigure[$C_1$,
  $V=-11.50$]{\includegraphics[width=3cm]{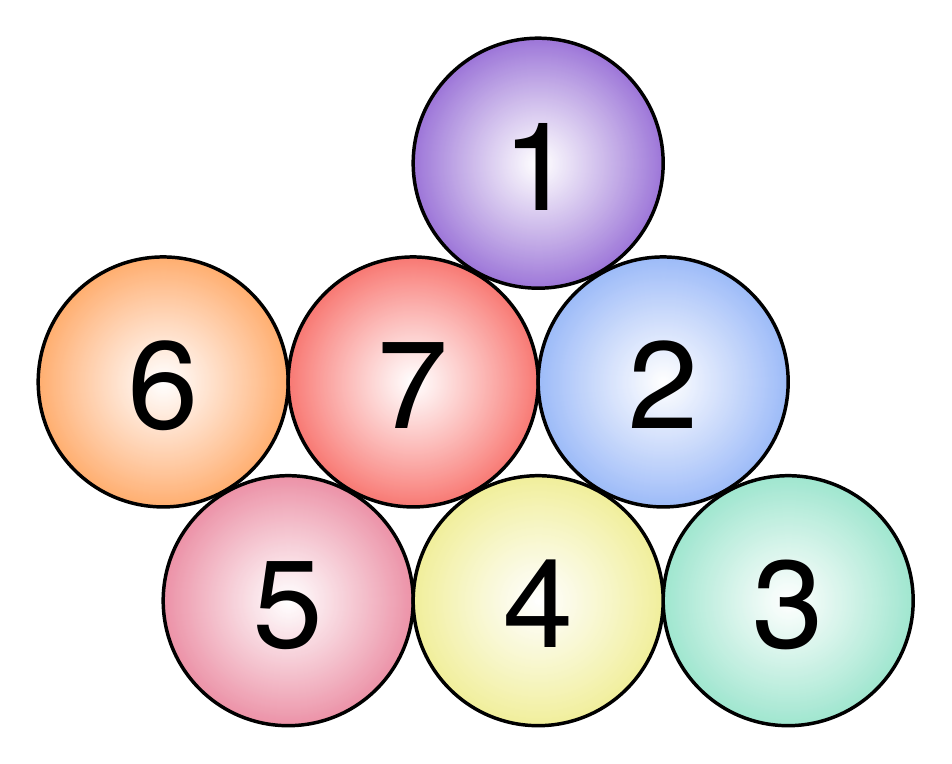}}
  \subfigure[$C_2$,
  $V=-11.48$]{\includegraphics[width=3cm]{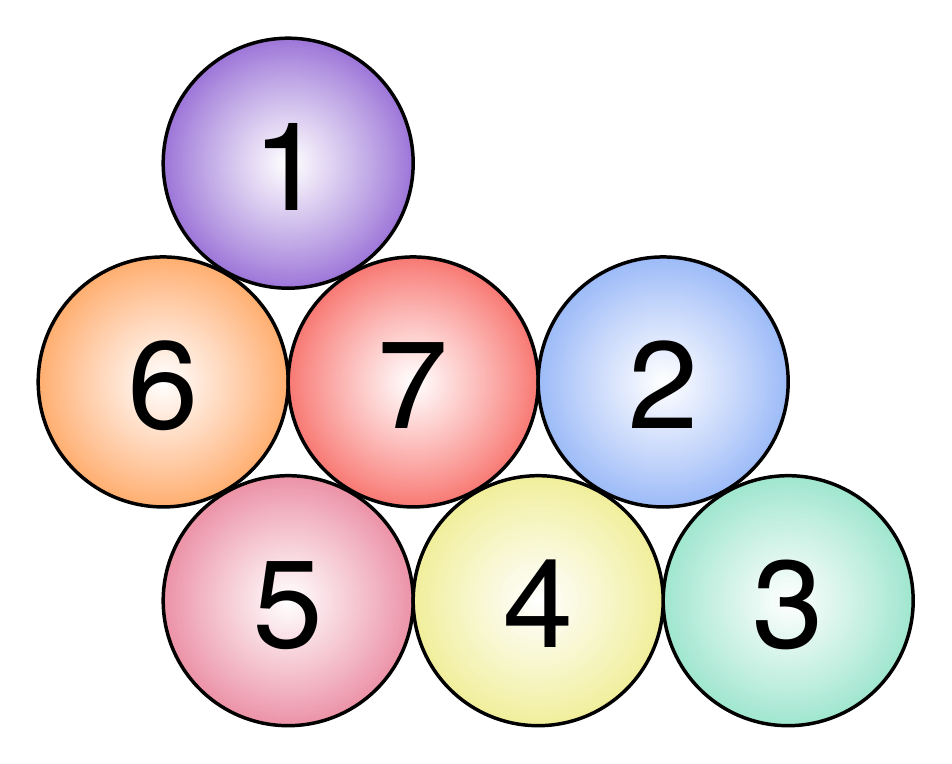}}
  \subfigure[$C_3$,
  $V=-11.40$]{\includegraphics[width=3cm]{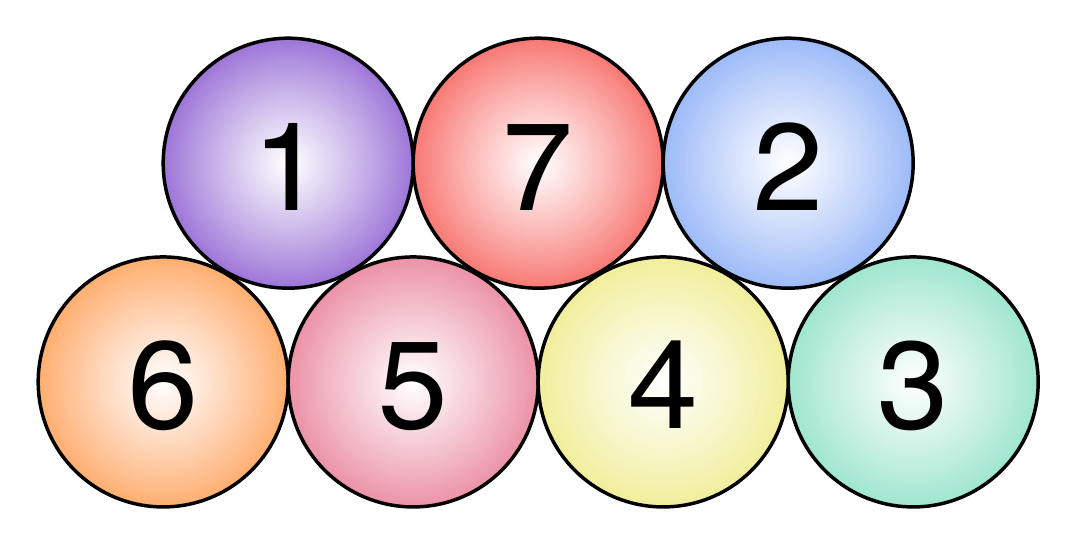}}
  \caption{Low energy conformations of the Lennard-Jones cluster of
    seven atoms in 2D.  Also indicated are the associated energy
    levels.} \label{f:lj7_conformations}
\end{figure}

The initial condition in this problem is near the closest packed
configuration, with six atoms located $\pi/3$ radians apart on a unit
circle and one in the center, as shown in Figure
\ref{f:lj7_conformations} (a).  We then look at first exits from this
configuration at inverse temperature $\beta =6$.  Exits correspond to
transitions into basins of attraction for a conformation other than
the lowest energy, closest packed, conformation; see
Figure~\ref{f:lj7_conformations}.  Here, basins of attraction are
those associated to the simple gradient dynamics $\dot{x}=-\nabla V
(x)$: these basins define the states on which ParRep is applied. For
the numbering of the conformations, we adopt the notation found in
\cite{Dellago:1998vs}, where the authors also explore transitions in
$\LJ$.  We also refer to~\cite{Hairer:2012wz} for an algorithm for
computing transitions within Lennard-Jones clusters.

Before proceeding to our results, we make note of several things that
are specific to this problem.  States are identified and compared as
follows:
\begin{enumerate}
\item Given the current configuration $X_t$, a gradient descent is run
  (following the dynamics $\dot {x} = -\nabla V (x)$ with $x_0 = X_t$)
  until the norm of the gradient becomes sufficiently small.  This is
  accomplished using the RK45 time discretization scheme, up to the
  time when the $\ell^2$ norm of the gradient, relative to the value
  of $\abs{V}$, is smaller than $10^{-5}$;
\item Once the local minima is found, in order to identify the
  conformation, the lengths of the twenty one ``bonds'' between pairs
  of atoms are computed.  Let $b_{ij}$ denote the distance between the
  atoms $i$ and $j$.
\item The current and previous states are compared using the $\ell^1$
  norm between the 21 bond lengths:
  \[
  \text{Distance to $C_0$} = \sum_{1 \le i<j \le 7} \abs{b_{ij} -
    b_{ij}^0}
  \]
  where $b_{ij}^0$ are the bond lengths of the $C_0$ conformation.
\item If this distance exceeds .5, the states are determined to be
  different.
\end{enumerate}
Next, when examining the output, we group conformations according to
energy levels.  In Figure \ref{f:lj7_conformations}, we present for
each conformation $C_0$, $C_1$, $C_2$ and $C_3$ a particular
numbering. Other permutations of the atoms also correspond to the same
energy level, and thus to the same conformation.

The observables that are used in this problem are:
\begin{subequations}
  \label{e:lj7_0bs}
  \begin{align}
    \text{Energy: }&V(x)\\
    \text{Square distance to the center of mass: } & \sum_{i=1}^7
    \norm{x_i -
      x_{\rm cm}}_{\ell^2}^2\\
    \text{Distance to $C_0$: } & \sum_{1 \le i<j \le 7} \abs{b_{ij} -
      b_{ij}^0}
  \end{align}
\end{subequations}
It is essential to use observables that are rotation and translation
invariant, as there is nothing to prevent the cluster from drifting or
rotating ({\it i.e.} no atom is pinned).

A related remark is that, to not introduce additional auxiliary
parameters, no confining potential is used.  Because of this, the
cluster does not always change into one of the other conformations.
In some cases, an atom can drift away from the cluster, and the
quenched configuration corresponds to the isolated atom and the five
remaning atoms.  In others, the quenched configuration corresponds to
$C_0$, but with two exterior atoms exchanged.  These scenarios
occurred in less than .1\% of the $10^5$ realizations of all
experiments, both ParRep and unaccelerated.

Over 99.9\% of the transitions from $C_0$ are to $C_1$ or $C_2$. A
transition from $C_0$ to $C_3$ has not been observed. Since there are
thus only two ways to leave $C_0$, we use the 95\% Clopper-Pearson
confidence intervals~\cite{Johnson:2005aa} (well adapted to binomial
random variables) to assess the quality of the exit point distribution
in our statistical tests.

\begin{figure}
  \subfigure{\includegraphics[width=7cm]{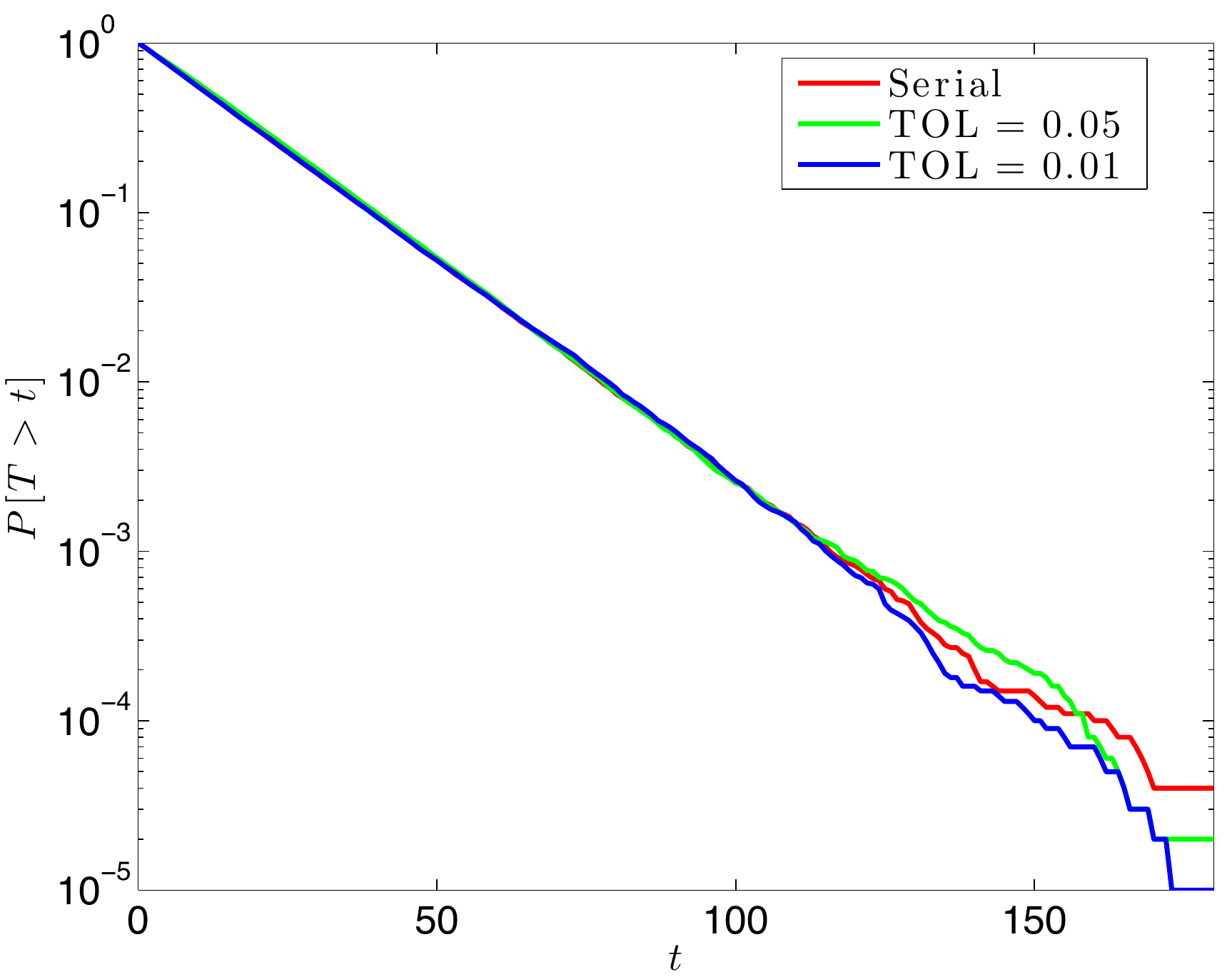}}
  \subfigure{\includegraphics[width=7cm]{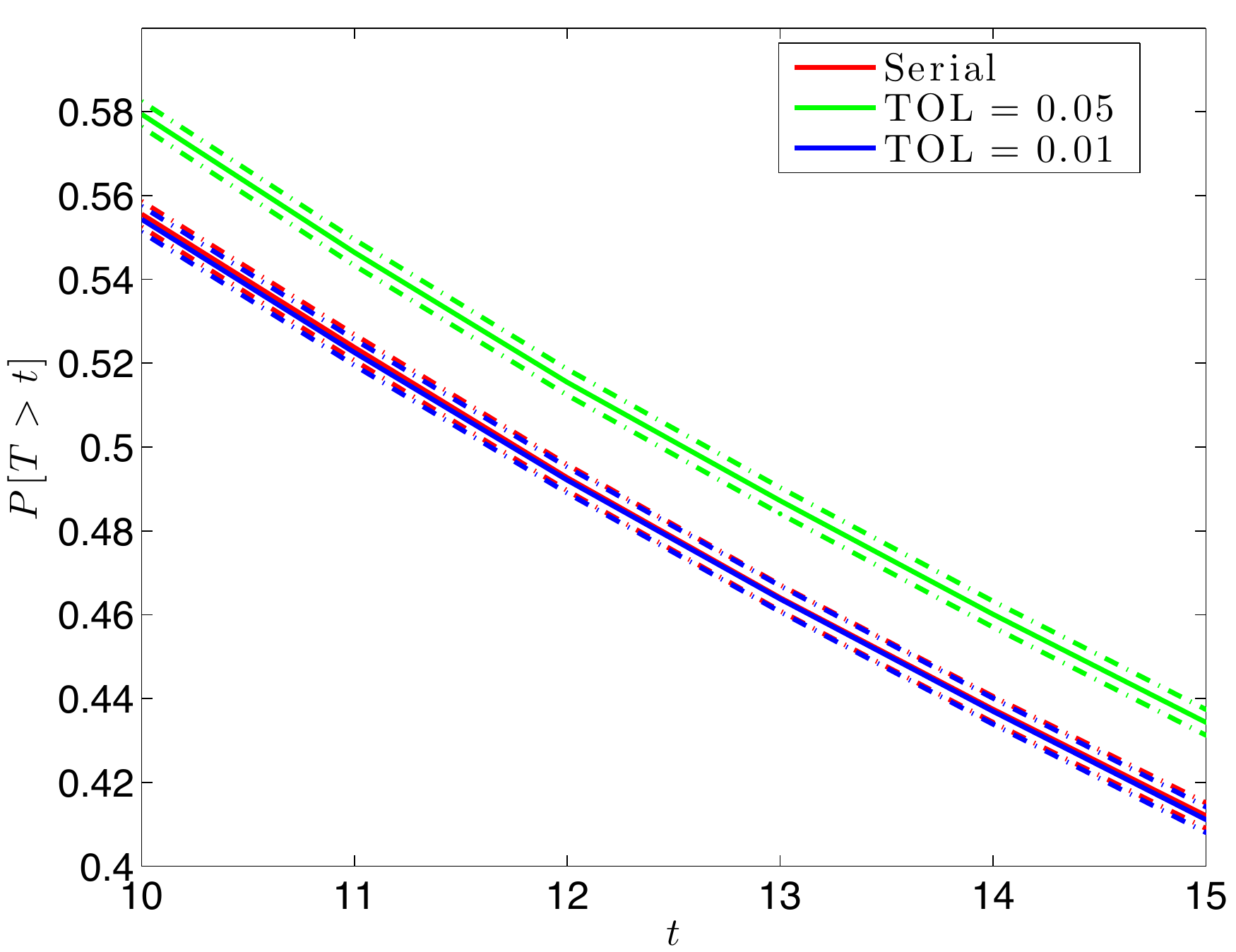}}
  \caption{{\bf $\LJ$-Single Escape:} Exit time distributions for
    $\LJ$ at $\beta = 6$.}
  \label{f:lj7beta6}
\end{figure}

Our statistical assessment of ParRep is given in Tables
\ref{t:lj7beta6_T} and \ref{t:lj7beta6_XT}.  It is only at the most
stringent tolerance of $\tol =0.01$ that an excellent agreement is
obtained in the exit time distribution, though some amount of speedup
is still gained. As in the other problems, the Kolmogorov-Smirnov test
may be an overly conservative measure of the quality of ParRep. As
shown in Figure \ref{f:lj7beta6}, the exit time distributions are
already in good qualitative agreement, even at less stringent
tolerances.  In contrast, the transition probabilities to the other
conformations are in very good agreement at tolerances beneath 0.1.

\begin{table}
  \caption{{\bf $\LJ$-Single Escape:} Comparison of the exit times for ParRep and an unaccelerated serial process
    escaping from the closest packed $\LJ$ configuration. }
  \label{t:lj7beta6_T}
  \begin{tabular}{l l c c c c c l}
    \hline \hline\\
    Method  & $\tol$ & $\mean{t_{\phase}}$ & $\Var(t_{\phase})$ & $\mean{T}$&
    $\mean{\text{Speedup}}$& \% Dephased & $T$ K.-S. Test ($p$)\\
    \hline\\
    Serial     & -- &--& --& 17.0 & --& --  &--\\
    ParRep & 0.2 & 0.411 & 0.0336 & 19.1 & 29.3  & 98.5\% &
    FAIL ($3.6 \times 10^{-235}$) \\
    ParRep &0.1 & .976& 0.125& 18.0 & 14.9 & 95.3\% & FAIL
    ($5.0 \times 10^{-55}$)\\
    ParRep & 0.05& 2.08 & 0.433& 17.6 & 7.83  & 90.0\% &
    FAIL ($9.4 \times
    10^{-28}$)\\
    ParRep &0.01 & 10.8 & 5.67&17.0 & 1.82 & 52.1\% &
    PASS (0.92)\\
    \hline
  \end{tabular}
\end{table}

\begin{table}
  \caption{{\bf $\LJ$-Single Escape:} Comparison of the hitting points for
    ParRep and an unaccelerated serial process
    escaping from the closest packed $\LJ$ configuration. }
  \label{t:lj7beta6_XT}
  \begin{tabular}{l l c c }
    \hline \hline\\
    Method  & $\tol$ &  $\prob[C_1]$  & $\prob[C_2]$ \\
    \hline\\
    Serial     & -- & (0.502, 0.508)& (0.491, 0.498) \\
    ParRep & 0.2 & (0.508, 0.514) & (0.485, 0.492) \\
    ParRep &0.1 & (0.506, 0.512)  & (0.488, 0.494) \\
    ParRep & 0.05& (0.505, 0.512) &  (0.488, 0.495)\\
    ParRep &0.01 & (0.504, 0.510) &  (0.490, 0.496) \\
    \hline
  \end{tabular}
\end{table}

\subsection{Conclusions from the Numerical Experiments}

From these numerical experiments, we observe that the modified ParRep
is indeed an efficient algorithm on various test cases, for which no a
priori knowledge on the decorrelation/dephasing time has been used.
{We thus trade setting $\tcorr$ and $\tphase$, {\it a
    priori}, for selecting physically informed observables together
  with a tolerance.  In conclusion, a tolerance on the order of 0.1
  in~\eqref{e:R2} seems to yield sensible results.  }

We also observe that the tolerance criteria used to assess
stationarity does not need to be very stringent to get the correct
distribution for the hitting points (thus to get the correct
state-to-state Markov chains).  In contrast, the K.-S. test on the
exit time distribution requires smaller tolerances to predict
statistical agreement.  However, even when the K.-S. test fails to
reject the null hypothesis, the agreement is very good.  Overall, it
is easier to obtain statistical agreement in the sequence of visited
states than it is to get statistical agreement of the time spent
within each state.


\bibliographystyle{plain}

\bibliography{acc_dyn}

\end{document}